\documentclass[a4paper]{amsart}
\usepackage{amssymb}
\usepackage[all]{xy}
\usepackage{enumitem}
\usepackage{nicefrac,stmaryrd}
\usepackage{mathtools}

\usepackage{graphicx}
\usepackage[pdftitle={Homological algebra in bivariant K-theory
  and other triangulated categories.  I},
  pdfauthor={Ralf Meyer and Ryszard Nest},
  pdfsubject={Mathematics; MSC 18E30, 19K35, 46L80, 55U35}
]{hyperref}
\usepackage[lite]{amsrefs}
\newcommand*{\MRref}[2]{\linebreak[0] \href{http://www.ams.org/mathscinet-getitem?mr=#1}{MR \textbf{#1}}}
\newcommand*{\arxiv}[1]{\linebreak[0] \href{http://www.arxiv.org/abs/#1}{arXiv:#1}}
\usepackage{microtype}

\DeclareMathOperator{\Hom}{Hom}

\DeclareMathOperator{\Res}{Res}

\DeclareMathOperator{\Tor}{Tor}
\DeclareMathOperator{\Ext}{Ext}

\DeclareMathOperator{\Rep}{Rep}

\DeclareMathOperator{\range}{range}
\DeclareMathOperator{\coker}{coker}

\DeclareMathOperator{\cone}{cone}

\newcommand*{\Left}{\mathbb L}
\newcommand*{\Right}{\mathbb R}
\newcommand*{\KK}{\textup{KK}}

\newcommand*{\K}{\textup K}

\newcommand*{\HH}{\textup{HH}}

\newcommand*{\Cst}{\textup C^*}

\newcommand*{\congto}{\xrightarrow\cong}
\newcommand*{\into}{\rightarrowtail}
\newcommand*{\prto}{\twoheadrightarrow}

\newcommand*{\Fam}{\mathcal F}% family of subgroups
\newcommand*{\CONT}{\textup C}% continuous functions

\newcommand*{\Cat}{\mathfrak C}% generic category, often Abelian
\newcommand*{\Tri}{\mathfrak T}% triangulated category
\newcommand*{\Ideal}{\mathfrak I}% ideal in a triangulated category
\newcommand*{\Proj}{\mathfrak P}% projective objects
\newcommand*{\Ab}{\mathfrak{Ab}}% category of Abelian groups
\newcommand*{\CAb}{\mathfrak{Ab}_\textup c}
\newcommand*{\FCAb}{\mathfrak{Ab}_\textup{fc}}
\newcommand*{\Mod}[1]{\mathfrak{Mod}(#1)}% category of modules
\newcommand*{\CMod}[1]{\mathfrak{Mod}(#1)_{\textup{c}}}
\newcommand*{\Fun}{\mathfrak{Fun}}% category of functors
\newcommand*{\Coh}{\mathfrak{Coh}}% category of finitely presented functors
\newcommand*{\Yoneda}{\mathbb Y}% Yoneda embedding
\newcommand*{\Ho}{\mathfrak{Ho}}% homotopy category
\newcommand*{\Der}{\mathfrak{Der}}% derived category

\newcommand*{\Hgy}{\textup H_*}
\newcommand*{\HGY}{\textup H}
\newcommand*{\IdealH}{\mathfrak I_{\HGY}}% kernel of the homology functor

\newcommand*{\VC}{{\mathcal{VC}}}

\newcommand*{\op}{\textup{op}}
\newcommand*{\red}{\textup{red}}
\newcommand*{\ID}{\textup{id}}

\newcommand*{\C}{\mathbb C}
\newcommand*{\R}{\mathbb R}
\newcommand*{\Sphere}{\mathbb S}
\newcommand*{\Torus}{\mathbb T}% torus
\newcommand*{\Z}{\mathbb Z}
\newcommand*{\Ztwo}{{\mathbb Z/2}}
\newcommand*{\N}{\mathbb N}
\newcommand*{\Comp}{\mathbb K}

\newcommand*{\nb}{\nobreakdash}

\mathtoolsset{mathic}

\DeclarePairedDelimiter{\sumclosure}{(}{)_\oplus}

\newcommand*{\lad}{\vdash}

\newcommand*{\rcross}{\mathbin{\ltimes_\textup r}\nobreak}
\newcommand*{\cross}{\mathbin\ltimes\nobreak}
\newcommand*{\blank}{\text\textvisiblespace}
\newcommand*{\inOb}{\mathrel{\in\in}\nobreak}

\newcommand*{\defeq}{\mathrel{\vcentcolon=}}

\numberwithin{equation}{section}

\theoremstyle{plain}
\newtheorem{theorem}[equation]{Theorem}
\newtheorem{proposition}[equation]{Proposition}
\newtheorem{lemma}[equation]{Lemma}

\newtheorem{corollary}[equation]{Corollary}
\theoremstyle{definition}
\newtheorem{definition}[equation]{Definition}
\newtheorem{notation}[equation]{Notation}
\theoremstyle{remark}
\newtheorem{remark}[equation]{Remark}
\newtheorem{warning}[equation]{Warning}
\newtheorem{example}[equation]{Example}

\begin{document}

\title[Homological algebra in bivariant \(\K\)-theory]{Homological algebra in
  bivariant \(\K\)-theory\\ and other triangulated categories.  I}

\author{Ralf Meyer}
\address{Mathematisches Institut\\
         Georg-August-Universit\"at G\"ottingen\\
         Bunsenstra\ss e 3--5\\
         37073 G\"ottingen\\
         Germany\\
}
\email{meyerr@member.ams.org}

\author{Ryszard Nest}
\address{K{\o}benhavns Universitets Institut for Matematiske Fag\\
         Universitetsparken 5\\ 2100 K{\o}benhavn\\ Denmark
}
\email{rnest@math.ku.dk}


\thanks{This research was supported by the EU-Network \emph{Quantum
    Spaces and Noncommutative Geometry} (Contract HPRN-CT-2002-00280).}

\begin{abstract}
  Bivariant (equivariant) \(\K\)\nb-theory is the standard setting for
  non-commutative topology.  We may carry over various techniques from
  homotopy theory and homological algebra to this setting.  Here we do
  this for some basic notions from homological algebra: phantom maps,
  exact chain complexes, projective resolutions, and derived functors.
  We introduce these notions and apply them to examples from bivariant
  \(\K\)\nb-theory.

  An important observation of Beligiannis is that we can approximate our
  category by an Abelian category in a canonical way, such that our
  homological concepts reduce to the corresponding ones in this Abelian
  category.  We compute this Abelian approximation in several
  interesting examples, where it turns out to be very concrete and
  tractable.

  The derived functors comprise the second page of a spectral
  sequence that, in favourable cases, converges towards
  Kasparov groups and other interesting objects.  This
  mechanism is the common basis for many different spectral
  sequences.  Here we only discuss the very simple
  \(1\)\nb-dimensional case, where the spectral sequences
  reduce to short exact sequences.
\end{abstract}

\maketitle

\section{Introduction}
\label{sec:intro}

It is well-known that many basic constructions from homotopy theory
extend to categories of \(\Cst\)\nb-algebras. As we argued
in~\cite{Meyer-Nest:BC}, the framework of \emph{triangulated categories}
is ideal for this purpose.  The notion of triangulated category was
introduced by Jean-Louis Verdier to formalise the properties of the
derived category of an Abelian category.  Stable homotopy theory
provides further classical examples of triangulated categories.  The
triangulated category structure encodes basic information about
manipulations with long exact sequences and (total) derived functors.
The main point of~\cite{Meyer-Nest:BC} is that the domain of the
Baum--Connes assembly map is the total left derived functor of the
functor that maps a \(G\)\nb-\(\Cst\)-algebra~\(A\) to \(\K_*(G\rcross
A)\).

The relevant triangulated categories in non-commutative topology come
from Kasparov's bivariant \(\K\)\nb-theory.  This bivariant version of
\(\K\)\nb-theory carries a composition product that turns it into a
category.  The formal properties of this and related categories are
surveyed in~\cite{Meyer:KK-survey}, with an audience of homotopy
theorists in mind.

\emph{Projective resolutions} are among the most fundamental concepts in
homological algebra; several others like derived functors are based on
it.  Projective resolutions seem to live in the underlying Abelian
category and not in its derived category.  This is why \emph{total}
derived functor make more sense in triangulated categories than the
derived functors themselves.  Nevertheless, we can define derived
functors in triangulated categories and far more general categories.
This goes back to Samuel Eilenberg and John~C. Moore
(\cite{Eilenberg-Moore:Foundations}).  We learned about this theory in
articles by Apostolos Beligiannis (\cite{Beligiannis:Relative}) and
J.~Daniel Christensen (\cite{Christensen:Ideals}).

Homological algebra in non-Abelian categories is always \emph{relative},
that is, we need additional structure to get started.  This is useful
because we may fit the additional data to our needs.  In a triangulated
category~\(\Tri\), there are several kinds of additional data that yield
equivalent theories; following~\cite{Christensen:Ideals}, we use an
\emph{ideal} in~\(\Tri\).  We only consider ideals~\(\Ideal\)
in~\(\Tri\) of the form
\[
\Ideal(A,B) \defeq \{x\in\Tri(A,B)\mid F(x)=0\}
\]
for a stable homological functor~\(F\colon \Tri\to\Cat\) into a stable
Abelian category~\(\Cat\).  Here \emph{stable} means that~\(\Cat\)
carries a suspension automorphism and that~\(F\) intertwines the
suspension automorphisms on \(\Tri\) and~\(\Cat\), and homological means
that exact triangles yield exact sequences.  Ideals of this form are
called \emph{homological ideals}.

A basic example is the ideal in the Kasparov category~\(\KK\) defined by
\begin{equation}
  \label{eq:def_Ideal_K}
  \Ideal_\K(A,B) \defeq
  \{f\in\KK(A,B)\mid 0=\K_*(f)\colon \K_*(A)\to\K_*(B)\}.
\end{equation}

For a compact (quantum) group~\(G\), we define two ideals
\(\Ideal_\ltimes\subseteq \Ideal_{\ltimes,\K} \subseteq \KK^G\) in the
equivariant Kasparov category~\(\KK^G\) by
\begin{align}
  \label{eq:def_Ideal_ltimes}
  \Ideal_\ltimes(A,B) &\defeq
  \{f\in\KK^G(A,B)\mid
  \text{\(G\cross f=0\) in \(\KK(G\cross A,G\cross B)\)}\},\\
  \label{eq:def_Ideal_ltimes_K}
  \Ideal_{\ltimes,\K}(A,B) &\defeq
  \{f\in\KK^G(A,B)\mid \K_*(G\cross f)=0\},
\end{align}
where \(\K_*(G\cross f)\) denotes the map~\(\K_*(G\cross A)\to
\K_*(G\cross B)\) induced by~\(f\).

For a locally compact group~\(G\) and a (suitable) family of
subgroups~\(\Fam\), we define the homological ideal
\begin{multline}
  \label{eq:def_VC}
  \VC_\Fam(A,B) \defeq \{f\in\KK^G(A,B)\mid \\
  \text{\(\Res_G^H(f)=0\) in \(\KK^H(A,B)\) for all \(H\in\Fam\)}\}.
\end{multline}
If~\(\Fam\) is the family of compact subgroups, then \(\VC_\Fam\) is
related to the Baum--Connes assembly map (\cite{Meyer-Nest:BC}).  Of
course, there are analogous ideals in more classical categories of
(spectra of) \(G\)\nb-CW-complexes.

All these examples can be analysed using the machinery we explain.  We
carry this out in some cases in Sections \ref{sec:plain_UCT}
and~\ref{sec:crossed_cqg}.

We use an ideal~\(\Ideal\) to carry over various notions from
homological algebra to our triangulated category~\(\Tri\).  In order
to see what they mean in examples, we characterise them using a stable
homological functor~\(F\colon \Tri\to\Cat\) with \(\ker F=\Ideal\).
This is often easy.  For instance, a chain complex with entries
in~\(\Tri\) is \emph{\(\Ideal\)\nb-exact} if and only if~\(F\) maps it
to an exact chain complex in the Abelian category~\(\Cat\) (see
Lemma~\ref{lem:exact_kerF}), and a morphism in~\(\Tri\) is an
\emph{\(\Ideal\)\nb-epimorphism} if and only if~\(F\) maps it to an
epimorphism.  Here we may take any functor~\(F\) with \(\ker
F=\Ideal\).

But the most crucial notions like projective objects and resolutions
require a more careful choice of the functor~\(F\).  Here we need the
\emph{universal \(\Ideal\)\nb-exact functor}, which is a stable
homological functor~\(F\) with \(\ker F=\Ideal\) such that any other
such functor factors uniquely through~\(F\) (up to natural equivalence).
The universal \(\Ideal\)\nb-exact functor and its applications are due
to Apostolos Beligiannis (\cite{Beligiannis:Relative}).

If \(F\colon \Tri\to\Cat\) is universal, then~\(F\) detects
\(\Ideal\)\nb-projective objects, and it identifies
\(\Ideal\)\nb-derived functors with derived functors in the Abelian
category~\(\Cat\) (see Theorem~\ref{the:universal_homological_nice}).
Thus all our homological notions reduce to their counterparts in the
Abelian category~\(\Cat\).

In order to apply this, we need to know when a functor~\(F\) with
\(\ker F=\Ideal\) is the universal one.  We develop a new, useful
criterion for this purpose here, which uses partially defined adjoint
functors (Theorem~\ref{the:universal_functor}).

Our criterion shows that the universal \(\Ideal_\K\)\nb-exact functor
for the ideal~\(\Ideal_\K\) in~\(\KK\) in~\eqref{eq:def_Ideal_K} is the
\(\K\)\nb-theory functor~\(\K_*\), considered as a functor from~\(\KK\)
to the category~\(\CAb^\Ztwo\) of countable \(\Ztwo\)\nb-graded Abelian
groups (see Theorem~\ref{the:K_projectives}).  Hence the derived
functors for~\(\Ideal_\K\) only involve \(\Ext\) and \(\Tor\) for
Abelian groups.

For the ideal~\(\Ideal_{\ltimes,\K}\) in~\(\KK^G\)
in~\eqref{eq:def_Ideal_ltimes_K}, we get the functor
\begin{equation}
  \label{eq:universal_Ideal_ltimes_K}
  \KK^G\to \CMod{\Rep G}^\Ztwo,\qquad
  A\mapsto \K_*(G\cross A),
\end{equation}
where \(\CMod{\Rep G}^\Ztwo\) denotes the Abelian category of
countable \(\Ztwo\)\nb-graded modules over the representation ring
\(\Rep G\) of the compact (quantum) group~\(G\) (see
Theorem~\ref{the:cross_K_universal}); here we use a certain canonical
\(\Rep G\)-module structure on \(\K_*(G\cross A)\).  Hence derived
functors with respect to~\(\Ideal_{\ltimes,\K}\) involve \(\Ext\) and
\(\Tor\) for \(\Rep G\)-modules.

We do not need the \(\Rep G\)-module structure on \(\K_*(G\cross A)\) to
define~\(\Ideal_{\ltimes,\K}\): our machinery notices automatically that
such a module structure is missing.  The universality of the functor
in~\eqref{eq:universal_Ideal_ltimes_K} clarifies in what sense
homological algebra with \(\Rep G\)-modules is a \emph{linearisation} of
algebraic topology with \(G\)\nb-\(\Cst\)-algebras.

The universal homological functor for the ideal~\(\Ideal_\ltimes\) is
quite similar to the one for~\(\Ideal_{\ltimes,\K}\) (see
Theorem~\ref{the:cross_universal}).  There is a canonical \(\Rep
G\)-module structure on \(G\cross A\) as an object of~\(\KK\), and the
universal \(\Ideal_\ltimes\)-exact functor is essentially the functor
\(A\mapsto G\cross A\), viewed as an object of a suitable Abelian
category that encodes this \(\Rep G\)-module structure; it also involves
a fully faithful embedding of~\(\KK\) in an Abelian category due to
Peter Freyd (\cite{Freyd:Representations}).

The derived functors that we have discussed above appear in a spectral
sequence which -- in favourable cases -- computes morphism spaces
in~\(\Tri\) (like \(\KK^G(A,B)\)) and other homological functors.  This
spectral sequence is a generalisation of the \emph{Adams spectral
  sequence} in stable homotopy theory and is the main motivation
for~\cite{Christensen:Ideals}.  Much earlier, such spectral sequences
were studied by Hans-Berndt Brinkmann in~\cite{Brinkmann:Relative}.
In~\cite{Meyer:Homology_in_KK_II}, this spectral sequence is applied to
our bivariant \(\K\)\nb-theory examples.  Here we only consider the much
easier case where this spectral sequence degenerates to an exact
sequence (see Theorem~\ref{the:UCT_homological}).  This generalises the
familiar Universal Coefficient Theorem for~\(\KK_*(A,B)\).

\section{Homological ideals in triangulated categories}
\label{sec:ideals_and_examples}

After fixing some basic notation, we introduce several interesting
ideals in bivariant Kasparov categories; we are going to discuss these
ideals throughout this article.  Then we briefly recall what a
triangulated category is and introduce homological ideals.  Before we
begin, we should point out that the choice of ideal is important because
all our homological notions depend on it.  It seems to be a matter of
experimentation and experience to find the right ideal for a given
purpose.

\subsection{Generalities about ideals in additive categories}
\label{sec:general_ideals}

All categories we consider will be \emph{additive}, that is, they have a
zero object and finite direct products and coproducts which agree, and
the morphism spaces carry Abelian group structures such that the
composition is additive in each variable (\cite{MacLane:Categories}).

\begin{notation}
  \label{note:morphisms}
  Let~\(\Cat\) be an additive category.  We write \(\Cat(A,B)\) for the
  group of morphisms \(A\to B\) in~\(\Cat\), and \(A\inOb\Cat\) to
  denote that~\(A\) is an object of the category~\(\Cat\).
\end{notation}

\begin{definition}
  \label{def:ideal}
  An \emph{ideal}~\(\Ideal\) in~\(\Cat\) is a family of subgroups
  \(\Ideal(A,B)\subseteq\Cat(A,B)\) for all \(A,B\inOb\Cat\) such that
  \[
  \Cat(C,D)\circ\Ideal(B,C)\circ\Cat(A,B)\subseteq\Ideal(A,D)
  \qquad\text{for all \(A,B,C,D\inOb \Cat\).}
  \]
\end{definition}

We write \(\Ideal_1\subseteq\Ideal_2\) if
\(\Ideal_1(A,B)\subseteq\Ideal_2(A,B)\) for all \(A,B\inOb \Cat\).
Clearly, the ideals in~\(\Tri\) form a complete lattice.  The largest
ideal~\(\Cat\) consists of all morphisms in~\(\Cat\); the smallest
ideal~\(0\) contains only zero morphisms.

\begin{definition}
  \label{def:kernel_ideal}
  Let \(\Cat\) and~\(\Cat'\) be additive categories and let \(F\colon
  \Cat\to\Cat'\) be an additive functor.  Its \emph{kernel}~\(\ker F\)
  is the ideal in~\(\Cat\) defined by
  \[
  \ker F(A,B) \defeq \{f\in\Cat(A,B) \mid F(f)=0\}.
  \]
\end{definition}

This should be distinguished from the \emph{kernel on objects},
consisting of all objects with \(F(A)\cong0\), which is used much more
frequently.  This agrees with the class of \(\ker F\)-contractible
objects that we introduce below.

\begin{definition}
  \label{def:quotient_category}
  Let \(\Ideal\subseteq\Tri\) be an ideal.  Its \emph{quotient
    category}~\(\Cat/\Ideal\) has the same objects as~\(\Cat\) and
  morphism groups \(\Cat(A,B)/\Ideal(A,B)\).
\end{definition}

The quotient category is again additive, and the obvious functor
\(F\colon \Cat\to\Cat/\Ideal\) is additive and satisfies \(\ker
F=\Ideal\).  Thus any ideal~\(\Ideal\) in~\(\Cat\) is of the form~\(\ker
F\) for a canonical additive functor~\(F\).

The additivity of \(\Cat/\Ideal\) and~\(F\) depends on the fact that any
ideal~\(\Ideal\) is compatible with \emph{finite} products in the
following sense: the natural isomorphisms
\begin{alignat*}{2}
  \Cat(A,B_1\times B_2) &\congto \Cat(A,B_1)\times\Cat(A,B_2),
  &\quad
  \Cat(A_1\times A_2,B) &\congto \Cat(A_1,B)\times\Cat(A_2,B)
  \\\intertext{restrict to isomorphisms}
  \Ideal(A,B_1\times B_2) &\congto \Ideal(A,B_1)\times\Ideal(A,B_2),
  &\quad
  \Ideal(A_1\times A_2,B) &\congto \Ideal(A_1,B)\times\Ideal(A_2,B).
\end{alignat*}

\subsection{Examples of ideals}
\label{sec:examples}

\begin{example}
  \label{exa:Kvan}
  Let~\(\KK\) be the \emph{Kasparov category}, whose objects are the
  separable \(\Cst\)\nb-algebras and whose morphism spaces are the
  Kasparov groups~\(\KK_0(A,B)\), with the Kasparov product as
  composition.  Let \(\Ab^\Ztwo\) be the category of
  \(\Ztwo\)\nb-graded Abelian groups.  Both categories are evidently
  additive.

  \(\K\)\nb-theory is an additive functor \(\K_*\colon
  \KK\to\Ab^\Ztwo\).  We let \(\Ideal_\K\defeq \ker \K_*\) (as
  in~\eqref{eq:def_Ideal_K}).  Thus \(\Ideal_\K(A,B)\subseteq\KK(A,B)\)
  is the kernel of the natural map
  \[
  \gamma\colon \KK(A,B) \to \Hom\bigl(\K_*(A),\K_*(B)\bigr)
  \defeq \prod_{n\in\Ztwo} \Hom\bigl(\K_n(A),\K_n(B)\bigr).
  \]

  There is another interesting ideal in~\(\KK\), namely, the kernel of a
  natural map
  \[
  \kappa \colon \Ideal_\K(A,B) \to \Ext\bigl(\K_*(A),\K_{*+1}(B)\bigr)
  \defeq \prod_{n\in\Ztwo} \Ext\bigl(\K_n(A),\K_{n+1}(B)\bigr)
  \]
  due to Lawrence Brown (see~\cite{Rosenberg-Schochet:UCT}), whose
  definition we now recall.  We represent \(f\in\KK(A,B)\cong
  \Ext\bigl(A, \CONT_0(\R,B)\bigr)\) by a \(\Cst\)\nb-algebra extension
  \(\CONT_0(\R,B)\otimes\Comp\into E\prto A\).  This yields an exact
  sequence
  \begin{equation}
    \label{eq:Ktheory_six-term}
    \begin{gathered}
      \xymatrix{
        \K_1(B) \ar[r] & \K_0(E) \ar[r] & \K_0(A) \ar[d]^{f_*} \\
        \K_1(A) \ar[u]_{f_*} & \K_1(E) \ar[l] & \K_0(B). \ar[l]
      }
    \end{gathered}
  \end{equation}
  The vertical maps in~\eqref{eq:Ktheory_six-term} are the two
  components of~\(\gamma(f)\).  If \(f\in\Ideal_\K(A,B)\),
  then~\eqref{eq:Ktheory_six-term} splits into two extensions of Abelian
  groups, which yield an element~\(\kappa(f)\) in
  \(\Ext\bigl(\K_*(A),\K_{*+1}(B)\bigr)\).
\end{example}

\begin{example}
  \label{exa:VC}
  Let~\(G\) be a second countable, locally compact group.  Let~\(\KK^G\)
  be the associated \emph{equivariant Kasparov category}; its objects
  are the separable \(G\)\nb-\(\Cst\)-algebras and its morphism spaces
  are the groups~\(\KK^G(A,B)\), with the Kasparov product as
  composition.  If \(H\subseteq\nobreak G\) is a closed subgroup, then
  there is a \emph{restriction functor} \(\Res_G^H\colon
  \KK^G\to\KK^H\), which simply forgets part of the equivariance.

  If~\(\Fam\) is a set of closed subgroups of~\(G\), we define an
  ideal~\(\VC_\Fam\) in~\(\KK^G\) by
  \[
  \VC_\Fam(A,B)\defeq \{f\in\KK^G(A,B)\mid \text{\(\Res_G^H(f)=0\) for
    all \(H\in\Fam\)}\}
  \]
  as in~\eqref{eq:def_VC}.  Of course, the condition \(\Res_G^H(f)=0\)
  is supposed to hold in~\(\KK^H(A,B)\).  We are mainly interested in
  the case where~\(\Fam\) is the family of all compact subgroups
  of~\(G\) and simply denote the ideal by~\(\VC\) in this case.

  This ideal arises if we try to compute \(G\)\nb-equivariant homology
  theories in terms of \(H\)\nb-equivariant homology theories for
  \(H\in\Fam\).  The ideal~\(\VC\) is closely related to the approach to
  the Baum--Connes assembly map in~\cite{Meyer-Nest:BC}.
\end{example}

The authors feel more at home with Kasparov theory than with spectra.
Many readers will prefer to work in categories of spectra of, say,
\(G\)\nb-CW-complexes.  We do not introduce these categories here; but
it shoud be clear enough that they support similar restriction functors,
which provide analogues of the ideals~\(\VC_\Fam\).

\begin{example}
  \label{exa:cross}
  Let \(G\) and \(\KK^G\) be as in Example~\ref{exa:VC}.  Using the
  crossed product functor (also called \emph{descent} functor)
  \[
  G\cross\blank\colon \KK^G\to\KK,\qquad
  A\mapsto G\cross A,
  \]
  we define ideals
  \(\Ideal_\ltimes\subseteq\Ideal_{\ltimes,\K}\subseteq\KK^G\) as in
  \eqref{eq:def_Ideal_ltimes} and~\eqref{eq:def_Ideal_ltimes_K} by
  \begin{align*}
    \Ideal_\ltimes(A,B) &\defeq \{f\in\KK^G(A,B)\mid
    \text{\(G\cross f=0\) in \(\KK(G\cross A,G\cross B)\)}\},\\
    \Ideal_{\ltimes,\K}(A,B) &\defeq \{f\in\KK^G(A,B)\mid
    \K_*(G\cross f)=0\colon \K_*(G\cross A)\to \K_*(G\cross B)\}.
  \end{align*}
  We only study these ideals for compact~\(G\).  In this case, the
  \emph{Green--Julg Theorem} identifies \(\K_*(G\cross A)\) with the
  \(G\)\nb-equivariant \(\K\)\nb-theory~\(\K_*^G(A)\)
  (see~\cite{Julg:K_equivariante}).  Hence the ideal
  \(\Ideal_{\ltimes,\K}\subseteq\KK^G\) is a good equivariant analogue
  of the ideal~\(\Ideal_\K\) in~\(\KK\).

  Literally the same definition as above provides ideals
  \(\Ideal_\ltimes\subseteq \Ideal_{\ltimes,\K}\subseteq \KK^G\)
  if~\(G\) is a compact \emph{quantum} group.  We will always allow this
  more general situation below, but readers unfamiliar with quantum
  groups may ignore this.
\end{example}

\begin{remark}
  \label{rem:quantum_group}%
  We emphasise quantum groups here because Examples \ref{exa:VC}
  and~\ref{exa:cross} become closely related in this context.  This
  requires a quantum group analogue of the ideals~\(\VC_\Fam\)
  in~\(\KK^G\) of Example~\ref{exa:VC}.  If~\(G\) is a locally compact
  quantum group, then Saad Baaj and Georges Skandalis construct a
  \(G\)\nb-equivariant Kasparov category~\(\KK^G\)
  in~\cite{Baaj-Skandalis:Hopf_KK}.  There is a forgetful functor
  \(\Res_G^H\colon \KK^G\to\KK^H\) for each closed quantum subgroup
  \(H\subseteq G\).  Therefore, a family~\(\Fam\) of closed quantum
  subgroups yields an ideal~\(\VC_\Fam\) in~\(\KK^G\) as in
  Example~\ref{exa:VC}.

  Let~\(G\) be a compact group as in Example~\ref{exa:cross}.  Any
  crossed product~\(G\cross A\) carries a canonical \emph{coaction}
  of~\(G\), that is, a coaction of the discrete quantum group
  \(\Cst(G)\).  \emph{Baaj-Skandalis duality} asserts that this yields an
  equivalence of categories \(\KK^G\cong \KK^{\Cst(G)}\)
  (see~\cite{Baaj-Skandalis:Hopf_KK}).  We get back the crossed product
  functor \(\KK^G\to\KK\) by composing this equivalence with the
  restriction functor \(\KK^{\Cst(G)}\to\KK\) for the trivial quantum
  subgroup.  Hence \(\Ideal_\ltimes\subseteq\KK^G\) corresponds by
  Baaj-Skandalis duality to \(\VC_\Fam\subseteq\KK^{\Cst(G)}\),
  where~\(\Fam\) consists only of the trivial quantum subgroup.

  Thus the constructions in Examples \ref{exa:VC} and~\ref{exa:cross}
  are both special cases of a more general construction for locally
  compact quantum groups.
\end{remark}

Finally, we consider a classic example from homological algebra.

\begin{example}
  \label{exa:homology}
  Let~\(\Cat\) be an Abelian category.  Let \(\Ho(A)\) be the homotopy
  category of unbounded chain complexes
  \[
  \dotsb \to
  C_n \xrightarrow{\delta_n}
  C_{n-1} \xrightarrow{\delta_{n-1}}
  C_{n-2} \xrightarrow{\delta_{n-2}}
  C_{n-3} \to \dotsb
  \]
  over~\(\Cat\).  The space of morphisms~\(A\to B\) in \(\Ho(\Cat)\) is
  the space~\([A,B]\) of \emph{homotopy classes} of chain maps
  from~\(A\) to~\(B\).

  Taking homology defines functors~\(H_n\colon \Ho(\Cat)\to\Cat\) for
  \(n\in\Z\), which we combine to a single functor \(\Hgy\colon
  \Ho(\Cat)\to \Cat^\Z\).  We let \(\IdealH\subseteq \Ho(\Cat)\) be its
  kernel:
  \begin{equation}
    \label{eq:kernel_homology}
    \IdealH(A,B)\defeq \{f\in [A,B]\mid \Hgy(f)=0\}.
  \end{equation}

  We also consider the category~\(\Ho(\Cat;\Z/p)\) of
  \emph{\(p\)\nb-periodic} chain complexes over~\(\Cat\) for
  \(p\in\N_{\ge1}\); its objects satisfy \(C_n=C_{n+p}\) and
  \(\delta_n=\delta_{n+p}\) for all \(n\in\Z\), and chain maps and
  homotopies are required \(p\)\nb-periodic as well.  The
  category~\(\Ho(\Cat;\Z/2)\) plays a role in connection with cyclic
  cohomology, especially with local cyclic cohomology
  (\cites{Puschnigg:Diffeotopy, Meyer:HLHA}).  The
  category~\(\Ho(\Cat;\Z/1)\) is isomorphic to the category of chain
  complexes without grading.  By convention, we let \(\Z/0=\Z\), so that
  \(\Ho(\Cat;\Z/0)=\Ho(\Cat)\).

  The homology of a periodic chain complex is, of course, periodic, so
  that we get a homological functor \(\Hgy\colon \Ho(\Cat;\Z/p)\to
  \Cat^{\Z/p}\); here~\(\Cat^{\Z/p}\) denotes the category of
  \(\Z/p\)\nb-graded objects of~\(\Cat\).  We let
  \(\IdealH\subseteq\Ho(\Cat;\Z/p)\) be the kernel of~\(\Hgy\) as
  in~\eqref{eq:kernel_homology}.
\end{example}

\subsection{What is a triangulated category?}
\label{sec:tri_cat}

A \emph{triangulated category} is a category~\(\Tri\) with a
\emph{suspension automorphism} \(\Sigma\colon \Tri\to\Tri\) and a class
of \emph{exact triangles}, subject to various axioms (see
\cites{Neeman:Triangulated, Meyer-Nest:BC, Verdier:Thesis}).  An exact
triangle is a diagram in~\(\Tri\) of the form
\[
A\to B\to C\to\Sigma A\qquad\text{or}\qquad
\xymatrix@C=1em@R=1.71em{
  A \ar[rr] & & B \ar[dl]\\
  &C, \ar[ul]|{[1]}
}
\]
where the~\([1]\) in the arrow~\(C\to A\) warns us that this map has
\emph{degree}~\(1\).  A \emph{morphism} of triangles is a triple of maps
\(\alpha,\beta,\gamma\) making the obvious diagram commute.

A typical example is the homotopy category~\(\Ho(\Cat;\Z/p))\) of
\(\Z/p\)\nb-graded chain complexes.  Here the suspension functor is
the (signed) \emph{translation} functor
\[
\begin{aligned}
  \Sigma\bigl((C_n,d_n)\bigr) &\defeq (C_{n-1},-d_{n-1})
  &\qquad&\text{on objects,}\\
  \Sigma\bigl((f_n)\bigr) &\defeq (f_{n-1})
  &\qquad&\text{on morphisms;}
\end{aligned}
\]
a triangle is exact if it is isomorphic to a \emph{mapping cone
  triangle}
\[
A\overset{f}\to B\to \cone(f) \to \Sigma A
\]
for some chain map~\(f\); the maps \(B\to\cone(f)\to\Sigma A\) are the
canonical ones.  It is well-known that this defines a triangulated
category for \(p=0\); the arguments for \(p\ge1\) are essentially the
same.

Another classical example is the stable homotopy category, say, of
compactly generated pointed topological spaces (it is not particularly
relevant which category of spaces or spectra we use).  The suspension is
\(\Sigma(A)\defeq \Sphere^1\wedge A\); a triangle is exact if it is
isomorphic to a \emph{mapping cone triangle}
\[
A\overset{f}\to B\to \cone(f) \to \Sigma A
\]
for some map~\(f\); the maps \(B\to\cone(f)\to\Sigma A\) are the
canonical ones.

We are mainly interested in the categories \(\KK\) and~\(\KK^G\)
introduced in~\S\ref{sec:examples}.  Their triangulated category
structure is discussed in detail in~\cite{Meyer-Nest:BC}.  We are facing
a notational problem because the functor \(X\mapsto \CONT_0(X)\) from
pointed compact spaces to \(\Cst\)\nb-algebras is \emph{contravariant},
so that \emph{mapping cone triangles} now have the form
\[
A \overset{f}\leftarrow B \leftarrow \cone(f) \leftarrow \CONT_0(\R,A)
\]
for a \(*\)\nb-homomorphism \(f\colon B\to A\); here
\[
\cone(f) = \bigl\{(a,b)\in \CONT_0\bigl((0,\infty],A\bigr)\times B\bigm|
a(\infty) = f(b)\bigr\}
\]
and the maps \(\CONT_0(\R,A)\to \cone(f)\to B\) are the obvious ones,
\(a\mapsto (a,0)\) and \((a,b)\mapsto b\).

It is reasonable to view a \(*\)\nb-homomorphism from~\(A\) to~\(B\)
as a morphism from~\(B\) to~\(A\).  Nevertheless, we prefer the
convention that an algebra homorphism \(A\to B\) is a morphism \(A\to
B\).  But then the most natural triangulated category structure lives on
the opposite category~\(\KK^\op\).  This creates only notational
difficulties because the opposite category of a triangulated category
inherits a canonical triangulated category structure, which has ``the
same'' exact triangles.  However, the passage to opposite categories
exchanges suspensions and desuspensions and modifies some sign
conventions.  Thus the functor \(A\mapsto \CONT_0(\R,A)\), which is the
suspension functor in~\(\KK^\op\), becomes the \emph{desuspension}
functor in~\(\KK\).  Fortunately, Bott periodicity implies that
\(\Sigma^2\cong\ID\), so that \(\Sigma\) and~\(\Sigma^{-1}\) agree.

Depending on your definition of a triangulated category, you may want
the suspension to be an equivalence or isomorphism of categories.  In
the latter case, you must replace~\(\KK^{(G)}\) by an equivalent
category (see~\cite{Meyer-Nest:BC}); since this is not important here,
we do not bother about this issue.

A triangle in~\(\KK^{(G)}\) is called \emph{exact} if it is isomorphic
to a mapping cone triangle
\[
\CONT_0(\R,B)\to \cone(f) \to A \overset{f}\to B
\]
for some (equivariant) \(*\)\nb-homomorphism~\(f\).

An important source of exact triangles in~\(\KK^G\) are
\emph{extensions}.  If \(A\into B\prto C\) is an extension of
\(G\)\nb-\(\Cst\)\nb-algebras with an equivariant completely positive
contractive section, then it yields a class in
\(\Ext(C,A)\cong\KK(\Sigma^{-1}C,A)\); the resulting triangle
\[
\Sigma^{-1}C\to A \to B\to C
\]
in~\(\KK^G\) is exact and called an \emph{extension triangle}.  It is
easy to see that any exact triangle is isomorphic to an extension
triangle.

It is shown in~\cite{Meyer-Nest:BC} that \(\KK\) and~\(\KK^G\) for a
locally compact group~\(G\) are triangulated categories with this extra
structure.  The same holds for the equivariant Kasparov theory~\(\KK^S\)
with respect to any \(\Cst\)\nb-bialgebra~\(S\); this theory was defined
by Baaj and Skandalis in~\cite{Baaj-Skandalis:Hopf_KK}.

The triangulated category axioms are discussed in greater detail in
\cites{Neeman:Triangulated, Meyer-Nest:BC, Verdier:Thesis}.  They encode
some standard machinery for manipulating long exact sequences.  Most of
them amount to formal properties of mapping cones and mapping cylinders,
which we can prove as in classical topology.  The only axiom that
requires more care is that any morphism \(f\colon A\to B\) should be
part of an exact triangle.

Unlike in~\cite{Meyer-Nest:BC}, we prefer to construct this triangle as
an extension triangle because this works in greater generality; we have
learned this idea from the work of Radu Popescu and Alexander Bonkat
(\cites{Popescu:Equivariant, Bonkat:Thesis}).  Any element in
\(\KK^S_0(A,B)\cong \KK^S_1\bigl(A,\CONT_0(\R,B)\bigr)\) can be
represented by an extension \(\Comp(\mathcal{H})\into E\prto A\) with an
equivariant completely positive contractive section, where
\(\mathcal{H}\) is a full \(S\)\nb-equivariant Hilbert
\(\CONT_0(\R,B)\)-module, so that \(\Comp(\mathcal{H})\) is
\(\KK^S\)-equivalent to \(\CONT_0(\R,B)\).  Hence the resulting
extension triangle in~\(\KK^S\) is isomorphic to one of the form
\[
\CONT_0(\R,A)\to \CONT_0(\R,B) \to E\to A;
\]
by construction, it contains the suspension of the given class in
\(\KK^S_0(A,B)\); it is easy to remove the suspension.

\begin{definition}
  \label{def:homological_functor}\label{def:cohomological_functor}
  Let~\(\Tri\) be a triangulated and~\(\Cat\) an Abelian category.  A
  covariant functor \(F\colon \Tri\to\Cat\) is called \emph{homological}
  if \(F(A)\to F(B)\to F(C)\) is exact at~\(F(B)\) for all exact
  triangles \(A\to B\to C\to \Sigma A\).  A contravariant functor with
  the analogous exactness property is called \emph{cohomological}.
\end{definition}

Let \(A\to B\to C\to \Sigma A\) be an exact triangle.  Then a homological
functor \(F\colon \Tri\to\Cat\) yields a natural long exact sequence
\[
\dotsb \to F_{n+1}(C)
\to F_n(A) \to F_n(B) \to F_n(C)
\to F_{n-1}(A) \to F_{n-1}(B)
\to \dotsb
\]
with \(F_n(A)\defeq F(\Sigma^{-n} A)\) for \(n\in\Z\), and a cohomological
functor \(F\colon \Tri^\op\to\Cat\) yields a natural long exact sequence
\[
\dotsb \leftarrow F^{n+1}(C)
\leftarrow F^n(A) \leftarrow F^n(B) \leftarrow F^n(C)
\leftarrow F^{n-1}(A) \leftarrow F^{n-1}(B)
\leftarrow \dotsb
\]
with \(F^n(A)\defeq F(\Sigma^{-n} A)\).

\begin{proposition}
  \label{pro:triangulated_rep_cohomological}
  Let~\(\Tri\) be a triangulated category.  The functors
  \[
  \Tri(A,\blank)\colon \Tri\to\Ab,
  \qquad
  B\mapsto \Tri(A,B)
  \]
  are homological for all \(A\inOb\Tri\).  Dually, the functors
  \[
  \Tri(\blank,B)\colon \Tri^\op\to\Ab,
  \qquad A\mapsto \Tri(A,B)
  \]
  are cohomological for all \(B\inOb\Tri\).
\end{proposition}

Observe that
\[
\Tri^n(A,B) = \Tri(\Sigma^{-n} A,B) \cong \Tri(A,\Sigma^{n}B) \cong
\Tri_{-n}(A,B).
\]

\begin{definition}
  \label{def:stable_homological}
  A \emph{stable additive category} is an additive category equipped
  with an (additive) automorphism~\(\Sigma\), called \emph{suspension}.

  A \emph{stable homological functor} is a homological functor \(F\colon
  \Tri\to\Cat\) into a stable Abelian category~\(\Cat\) together with
  natural isomorphisms \(F\bigl(\Sigma_\Tri(A)\bigr)\cong
  \Sigma_\Cat\bigl(F(A)\bigr)\) for all \(A\inOb\Tri\).
\end{definition}

\begin{example}
  \label{exa:homological_functors}
  The category~\(\Cat^{\Z/p}\) of \(\Z/p\)-graded objects of an
  Abelian category~\(\Cat\) is stable for any \(p\in\N\); the suspension
  automorphism merely shifts the grading.  The functors \(\K_*\colon
  \KK\to\Ab^\Ztwo\) and \(\Hgy\colon \Ho(\Cat;\Z/p)\to \Cat^{\Z/p}\)
  introduced in Examples \ref{exa:Kvan} and~\ref{exa:homology} are
  stable homological functors.
\end{example}

If \(F\colon \Tri\to\Cat\) is any homological functor, then
\[
F_*\colon \Tri\to\Cat^\Z,
\qquad A\mapsto \bigl(F_n(A)\bigr)_{n\in\Z}
\]
is a stable homological functor.  Many of our examples satisfy
\emph{Bott periodicity}, that is, there is a natural isomorphism
\(F_2(A)\cong F(A)\).  Then we get a stable homological functor
\(F_*\colon \Tri\to\Cat^{\Ztwo}\).  A typical example for this is the
functor~\(\K_*\).

\begin{definition}
  \label{def:exact_functor}
  A functor \(F\colon \Tri\to\Tri'\) between two triangulated categories
  is called \emph{exact} if it intertwines the suspension automorphisms
  (up to specified natural isomorphisms) and maps exact triangles
  in~\(\Tri\) again to exact triangles in~\(\Tri'\).
\end{definition}

\begin{example}
  \label{exa:exact_functors}
  The restriction functor \(\Res^H_G\colon \KK^G\to\KK^H\) for a closed
  quantum subgroup~\(H\) of a locally compact quantum group~\(G\) and
  the crossed product functors \(G\cross\blank, G\rcross\blank\colon
  \KK^G\to\KK\) are exact because they preserve mapping cone triangles.
\end{example}

Let \(F\colon \Tri_1\to\Tri_2\) be an exact functor.  If \(G\colon
\Tri_2\to?\) is exact, homological, or cohomological, then so is
\(G\circ F\).

Using Examples \ref{exa:homological_functors}
and~\ref{exa:exact_functors}, we see that the functors that define the
ideals \(\ker\gamma\) in Example~\ref{exa:Kvan}, \(\VC_\Fam\) in
Example~\ref{exa:VC}, \(\Ideal_\ltimes\), and~\(\Ideal_{\ltimes,\K}\) in
Example~\ref{exa:cross}, and~\(\IdealH\) in Example~\ref{exa:homology}
are all stable and either homological or exact.

\subsection{The universal homological functor}
\label{sec:universal_homological}

The following general construction of Peter Freyd
(\cite{Freyd:Representations}) plays an important role
in~\cite{Beligiannis:Relative}.  For an additive category~\(\Cat\), let
\(\Fun(\Cat^\op,\Ab)\) be the category of contravariant additive
functors \(\Cat\to\Ab\), with natural transformations as morphisms.
Unless~\(\Cat\) is essentially small, this is not quite a category
because the morphisms may form classes instead of sets.  We may ignore
this set-theoretic problem because the bivariant Kasparov categories
that we are interested in are essentially small, and the
subcategory~\(\Coh(\Cat)\) of \(\Fun(\Cat^\op,\Ab)\) that we are going
to use later on is an honest category for any~\(\Cat\).

The category \(\Fun(\Cat^\op,\Ab)\) is Abelian: if \(f\colon F_1\to
F_2\) is a natural transformation, then its kernel, cokernel, image, and
co-image are computed pointwise on the objects of~\(\Cat\), so that they
boil down to the corresponding constructions with Abelian groups.

The \emph{Yoneda embedding} is an additive functor
\[
\Yoneda\colon \Cat\to \Fun(\Cat^\op,\Ab),
\qquad
B\mapsto \Tri(\blank,B).
\]
This functor is fully faithful, and there are natural isomorphisms
\[
\Hom(\Yoneda(B),F)\cong F(B)
\qquad \text{for all \(F\inOb\Fun(\Cat^\op,\Ab)\), \(B\inOb\Tri\)}
\]
by the \emph{Yoneda lemma}.  A functor \(F\inOb\Fun(\Cat^\op,\Ab)\) is
called \emph{representable} if it is isomorphic to \(\Yoneda(B)\) for
some \(B\inOb\Cat\).  Hence~\(\Yoneda\) yields an equivalence of
categories between~\(\Cat\) and the subcategory of representable
functors in \(\Fun(\Cat^\op,\Ab)\).

A functor \(F\inOb\Fun(\Cat^\op,\Ab)\) is called \emph{finitely
  presented} if there is an exact sequence \(\Yoneda(B_1)\to
\Yoneda(B_2) \to F\to 0\) with \(B_1,B_2\inOb\Tri\).  Since~\(\Yoneda\)
is fully faithful, this means that~\(F\) is the cokernel of
\(\Yoneda(f)\) for a morphism~\(f\) in~\(\Cat\).  We let \(\Coh(\Cat)\)
be the full subcategory of finitely presented functors in
\(\Fun(\Cat^\op,\Ab)\).  Since representable functors belong to
\(\Coh(\Cat)\), we still have a Yoneda embedding \(\Yoneda\colon
\Cat\to\Coh(\Cat)\).  Although the category \(\Coh(\Tri)\) tends to be
very big and therefore unwieldy, it plays an important theoretical role.

\begin{theorem}[Freyd's Theorem]
  \label{the:Freyd}
  Let~\(\Tri\) be a triangulated category.

  Then \(\Coh(\Tri)\) is a stable Abelian category that has enough
  projective and enough injective objects, and the projective and
  injective objects coincide.

  The functor \(\Yoneda\colon \Tri\to \Coh(\Tri)\) is fully faithful,
  stable, and homological.  Its essential range \(\Yoneda(\Tri)\)
  consists of projective-injective objects.  Conversely, an object of
  \(\Coh(\Tri)\) is projective-injective if and only if it is a retract
  of an object of \(\Yoneda(\Tri)\).

  The functor~\(\Yoneda\) is the universal \textup{(}stable\textup{)}
  homological functor in the following sense: any
  \textup{(}stable\textup{)} homological functor \(F\colon
  \Tri\to\Cat'\) to a \textup{(}stable\textup{)} Abelian
  category~\(\Cat'\) factors uniquely as \(F=\bar{F}\circ\Yoneda\) for a
  \textup{(}stable\textup{)} exact functor \(F\colon
  \Coh(\Tri)\to\Cat'\).
\end{theorem}

If idempotents in~\(\Tri\) split -- as in all our examples -- then
\(\Yoneda(\Tri)\) is closed under retracts, so that \(\Yoneda(\Tri)\) is
equal to the class of projective-injective objects in \(\Coh(\Tri)\).

\subsection{Homological ideals in triangulated categories}
\label{sec:tri_ideals}

Let~\(\Tri\) be a triangulated category, let~\(\Cat\) be a stable
additive category, and let \(F\colon \Tri\to\Cat\) be a stable
homological functor.  Then \(\ker F\) is a stable ideal in the following
sense:

\begin{definition}
  \label{def:stable_ideal}
  An ideal~\(\Ideal\) in~\(\Tri\) is called \emph{stable} if the
  suspension isomorphisms \(\Sigma\colon \Tri(A,B)\congto\Tri(\Sigma
  A,\Sigma B)\) for \(A,B\inOb\Tri\) restrict to isomorphisms
  \[
  \Sigma\colon \Ideal(A,B)\congto\Ideal(\Sigma A,\Sigma B).
  \]
\end{definition}

If~\(\Ideal\) is stable, then there is a unique suspension automorphism
on~\(\Tri/\Ideal\) for which the canonical functor
\(\Tri\to\Tri/\Ideal\) is stable.  Thus the stable ideals are exactly
the kernels of stable additive functors.

\begin{definition}
  \label{def:homological_ideal}
  An ideal~\(\Ideal\) in a triangulated category~\(\Tri\) is called
  \emph{homological} if it is the kernel of a stable homological
  functor.
\end{definition}

\begin{remark}
  \label{rem:homological_without_exact_triangles}
  Freyd's Theorem shows that~\(\Yoneda\) induces a bijection between
  (stable) exact functors \(\Coh(\Tri)\to\Cat'\) and (stable)
  homological functors \(\Tri\to\Cat'\) because \(\bar{F}\circ\Yoneda\)
  is homological if \(\bar{F}\colon \Coh(\Tri)\to\Cat'\) is exact.
  Hence the notion of homological functor is independent of the
  triangulated category structure on~\(\Tri\) because the Yoneda
  embedding \(\Yoneda\colon \Tri\to\Coh(\Tri)\) does not involve any
  additional structure.  Hence the notion of homological ideal only uses
  the suspension automorphism, not the class of exact triangles.
\end{remark}

All the ideals considered in~\S\ref{sec:examples} except for
\(\ker\kappa\) in Example~\ref{exa:Kvan} are kernels of stable
homological functors or exact functors.  Those of the first kind are
homological by definition.  If \(F\colon\Tri\to \Tri'\) is an exact
functor between two triangulated categories, then \(\Yoneda\circ F\colon
\Tri\to\Coh(\Tri')\) is a stable homological functor with \(\ker
\Yoneda\circ F=\ker F\) by Freyd's Theorem~\ref{the:Freyd}.  Hence
kernels of exact functors are homological as well.

Is any homological ideal the kernel of an exact functor?  This is
\emph{not} the case:

\begin{proposition}
  \label{pro:counter_kernel_exact}
  Let \(\Der(\Ab)\) be the derived category of the category~\(\Ab\) of
  Abelian groups.  Define the ideal~\(\IdealH\) in \(\Der(\Ab)\) as in
  Example~\ref{exa:homology}.  This ideal is not the kernel of an exact
  functor.
\end{proposition}

We postpone the proof to the end
of~\S\ref{sec:basic_algebra_from_ideals} because it uses the machinery
of~\S\ref{sec:basic_algebra_from_ideals}.

It takes some effort to characterise homological ideals because
\(\Tri/\Ideal\) is almost never Abelian.  The results in
\cite{Beligiannis:Relative}*{\S2--3} show that an ideal is
homological if and only if it is \emph{saturated} in the notation
of~\cite{Beligiannis:Relative}.  We do not discuss this notion here
because most ideals that we consider are obviously homological.  The
only example where we could profit from an abstract characterisation is
the ideal~\(\ker\kappa\) in Example~\ref{exa:Kvan}.

There is no obvious homological functor whose kernel is~\(\ker\kappa\)
because~\(\kappa\) is not a functor on~\(\KK\).  Nevertheless,
\(\ker\kappa\) is the kernel of an exact functor; the relevant functor
is the functor \(\KK\to\operatorname{UCT}\), where
\(\operatorname{UCT}\) is the variant of~\(\KK\) that satisfies the
Universal Coefficient Theorem in complete generality.  This functor can
be constructed as a localisation of~\(\KK\) (see~\cite{Meyer-Nest:BC}).
The Universal Coefficient Theorem implies that its kernel is
exactly~\(\ker\kappa\).

\section{From homological ideals to derived functors}
\label{sec:homological_ideals_to_derived_functors}

Once we have a stable homological functor \(F\colon \Tri\to\Cat\), it is
not surprising that we can do a certain amount of homological algebra
in~\(\Tri\).  For instance, we may call a chain complex of objects
of~\(\Tri\) \emph{\(F\)\nb-exact} if~\(F\) maps it to an exact chain
complex in~\(\Cat\); and we may call an object
\emph{\(F\)\nb-projective} if~\(F\) maps it to a projective object
in~\(\Cat\).  But are these definitions reasonable?

We propose that a reasonable homological notion should depend only on
the ideal~\(\ker F\).  We will see that the notion of \(F\)\nb-exact
chain complex is reasonable and only depends on \(\ker F\).  In
contrast, the notion of projectivity above depends on~\(F\) and is only
reasonable in special cases.  There is another, more useful, notion of
projective object that depends only on the ideal~\(\ker F\).

Various notions from homological algebra still make sense in the context
of homological ideals in triangulated categories.  Our discussion mostly
follows \cites{Asadollahi-Salarian:Gorenstein, Beligiannis:Relative,
  Christensen:Ideals, Eilenberg-Moore:Foundations}.  All our definitions
involve only the ideal, not a stable homological functor that defines
it.  We reformulate them in terms of an exact or a stable homological
functor defining the ideal in order to understand what they mean in
concrete cases.  Following~\cite{Eilenberg-Moore:Foundations}, we
construct projective objects using adjoint functors.

The most sophisticated concept in this section is the \emph{universal
  \(\Ideal\)\nb-exact functor}, which gives us complete control over
projective resolutions and derived functors.  We can usually describe
such functors very concretely.

\subsection{Basic notions}
\label{sec:basic_algebra_from_ideals}

We introduce some useful terminology related to an ideal:

\begin{definition}
  \label{def:ideal_monic}\label{def:ideal_epic}\label{def:ideal_phantom_map}\label{def:ideal_equivalence}\label{def:ideal_contractible_object}\label{def:ideal_exact_triangle}
  Let~\(\Ideal\) be a homological ideal in a triangulated
  category~\(\Tri\).

  \begin{itemize}
  \item Let \(f\colon A\to B\) be a morphism in~\(\Tri\); embed it in an
    exact triangle \(A\overset{f}\to B\overset{g}\to C\overset{h}\to
    \Sigma A\).  We call~\(f\)
    \begin{itemize}
    \item \emph{\(\Ideal\)\nb-monic} if \(h\in\Ideal\);

    \item \emph{\(\Ideal\)\nb-epic} if \(g\in\Ideal\);

    \item an \emph{\(\Ideal\)\nb-equivalence} if it is both
      \(\Ideal\)\nb-monic and \(\Ideal\)\nb-epic, that is,
      \(g,h\in\Ideal\);

    \item an \emph{\(\Ideal\)\nb-phantom map} if \(f\in\Ideal\).
    \end{itemize}

  \item An object \(A\inOb\Tri\) is called
    \emph{\(\Ideal\)\nb-contractible} if \(\ID_A\in\Ideal(A,A)\).

  \item An exact triangle \(A\overset{f}\to B\overset{g}\to
    C\overset{h}\to \Sigma A\) in~\(\Tri\) is called
    \emph{\(\Ideal\)\nb-exact} if \(h\in\Ideal\).

  \end{itemize}
\end{definition}

The notions of monomorphism (or monic morphism) and epimorphism (or epic
morphism) -- which can be found in any book on category theory
such as~\cite{MacLane:Categories} -- are categorical ways to
express injectivity or surjectivity of maps.  A morphism in an Abelian
category that is both monic and epic is invertible.

The classes of \(\Ideal\)\nb-phantom maps, \(\Ideal\)\nb-monics,
\(\Ideal\)\nb-epics, and of \(\Ideal\)\nb-exact triangles determine
each other uniquely because we can embed any morphism in an exact
triangle in any position.  It is a matter of taste which of these is
considered most fundamental.  Following Daniel Christensen
(\cite{Christensen:Ideals}), we favour the phantom maps.  Other authors
prefer exact triangles instead (\cites{Asadollahi-Salarian:Gorenstein,
  Beligiannis:Relative, Eilenberg-Moore:Foundations}).  Of course, the
notion of an \(\Ideal\)\nb-phantom map is redundant; it becomes more
relevant if we consider, say, the class of \(\Ideal\)\nb-exact
triangles as our basic notion.

Notice that~\(f\) is \(\Ideal\)\nb-epic or \(\Ideal\)\nb-monic if
and only if~\(-f\) is.  If~\(f\) is \(\Ideal\)\nb-epic or
\(\Ideal\)\nb-monic, then so are \(\Sigma^n(f)\) for all \(n\in\Z\)
because~\(\Ideal\) is stable.  Similarly, (signed) suspensions of
\(\Ideal\)\nb-exact triangles remain \(\Ideal\)\nb-exact triangles.

\begin{lemma}
  \label{lem:monic_kerF}\label{lem:epic_kerF}\label{lem:ideal_exact_triangle_kerF}
  Let \(F\colon \Tri\to\Cat\) be a stable homological functor into a
  stable Abelian category~\(\Cat\).
  \begin{itemize}
  \item A morphism~\(f\) in~\(\Tri\) is
    \begin{itemize}
    \item a \(\ker F\)-phantom map if and only if \(F(f)=0\);

    \item \(\ker F\)-monic if and only if \(F(f)\) is monic;

    \item \(\ker F\)-epic if and only if \(F(f)\) is epic;

    \item a \(\ker F\)-equivalence if and only if \(F(f)\) is
      invertible.
    \end{itemize}

  \item An object \(A\inOb\Tri\) is \(\ker F\)-contractible if and
    only if \(F(A)=0\).

  \item An exact triangle \(A\to B\to C\to\Sigma A\) is \(\ker
    F\)-exact if and only if
    \[
    0 \to F(A)\to F(B)\to F(C)\to 0
    \]
    is a short exact sequence in~\(\Cat\).

  \end{itemize}
\end{lemma}

\begin{proof}
  Sequences in~\(\Cat\) of the form \(X\overset{0}\to Y\overset{f}\to
  Z\) or \(X\overset{f}\to Y\overset{0}\to Z\) are exact at~\(Y\) if and
  only if~\(f\) is monic or epic, respectively.  Moreover, a sequence of
  the form \(X\overset{0} \to Y\to Z \to U\overset{0}\to W\) is exact if
  and only if \(0\to Y\to Z\to U\to 0\) is exact.

  Combined with the long exact homology sequences for~\(F\) and suitable
  exact triangles, these observations yield the assertions about
  monomorphisms, epimorphisms, and exact triangles.  The description of
  equivalences and contractible objects follows, and phantom maps are
  trivial, anyway.
\end{proof}

Now we specialise these notions to the ideal \(\Ideal_\K\subseteq\KK\)
of Example~\ref{exa:Kvan}, replacing~\(\Ideal_\K\) by~\(\K\) in our
notation to avoid clutter.

\begin{itemize}
\item Let \(f\in\KK(A,B)\) and let \(\K_*(f)\colon \K_*(A)\to\K_*(B)\)
  be the induced map.  Then~\(f\) is
  \begin{itemize}
  \item a \emph{\(\K\)\nb-phantom map} if and only if \(\K_*(f)=0\);

  \item \emph{\(\K\)\nb-monic} if and only if \(\K_*(f)\) is
    injective;

  \item \emph{\(\K\)\nb-epic} if and only if \(\K_*(f)\) is
    surjective;

  \item a \emph{\(\K\)\nb-equivalence} if and only if \(\K_*(f)\) is
    invertible.
  \end{itemize}

\item A \(\Cst\)\nb-algebra \(A\inOb\KK\) is
  \emph{\(\K\)\nb-contractible} if and only if \(\K_*(A)=0\).

\item An exact triangle \(A\to B\to C\to \Sigma A\) in \(\KK\) is
  \emph{\(\K\)\nb-exact} if and only if
  \[
  0 \to \K_*(A) \to \K_*(B) \to \K_*(C) \to 0
  \]
  is a short exact sequence (of \(\Ztwo\)\nb-graded Abelian groups).
\end{itemize}

Similar things happen for the other ideals in~\S\ref{sec:examples} that
are \emph{naturally} defined as kernels of stable homological functors.

\begin{remark}
  \label{rem:warn_stable_homological_needed}
  It is crucial for the above theory that we consider functors that are
  both \emph{stable} and \emph{homological}.  Everything fails if we
  drop either assumption and consider functors such as \(\K_0(A)\) or
  \(\Hom\bigl(\Z/4,\K_*(A)\bigr)\).
\end{remark}

\begin{lemma}
  \label{lem:equivalence_contractible}
  An object \(A\inOb\Tri\) is \(\Ideal\)\nb-contractible if and only
  if \(0\colon 0\to A\) is an \(\Ideal\)\nb-equivalence.  A
  morphism~\(f\) in~\(\Tri\) is an \(\Ideal\)\nb-equivalence if and
  only if its generalised mapping cone is \(\Ideal\)\nb-contractible.
\end{lemma}

Thus the classes of \(\Ideal\)\nb-equivalences and of
\(\Ideal\)\nb-contractible objects determine each other.  But they do
not allow us to recover the ideal itself.  For instance, the ideals
\(\Ideal_\K\) and~\(\ker\kappa\) in Example~\ref{exa:Kvan} have the same
contractible objects and equivalences.

\begin{proof}
  Recall that the generalised mapping cone of~\(f\) is the object~\(C\)
  that fits in an exact triangle \(A\overset{f}\to B\to C\to \Sigma A\).
  The long exact sequence for this triangle yields that \(F(f)\) is
  invertible if and only if \(F(C)=0\), where~\(F\) is some stable
  homological functor~\(F\) with \(\ker F=\Ideal\).  Now the second
  assertion follows from Lemma~\ref{lem:epic_kerF}.  Since the
  generalised mapping cone of \(0\to A\) is~\(A\), the first assertion
  is a special case of the second one.
\end{proof}

Many ideals are defined as \(\ker F\) for an exact functor \(F\colon
\Tri\to\Tri'\) between triangulated categories.  We can also use such a
functor to describe the above notions:

\begin{lemma}
  \label{lem:monic_kerF_exact}\label{lem:epic_kerF_exact}\label{lem:ideal_exact_triangle_kerF_exact}
  Let \(\Tri\) and~\(\Tri'\) be triangulated categories and let
  \(F\colon \Tri\to\Tri'\) be an exact functor.
  \begin{itemize}
  \item A morphism \(f\in\Tri(A,B)\) is
    \begin{itemize}
    \item a \(\ker F\)-phantom map if and only if \(F(f)=0\);

    \item \(\ker F\)-monic if and only if \(F(f)\) is (split)
      monic.

    \item \(\ker F\)-epic if and only if \(F(f)\) is (split) epic;

    \item a \(\ker F\)-equivalence if and only if \(F(f)\) is
      invertible.
    \end{itemize}

  \item An object \(A\inOb\Tri\) is \(\ker F\)-contractible if and
    only if \(F(A)=0\).

  \item An exact triangle \(A\to B\to C\to\Sigma A\) is \(\ker
    F\)-exact if and only if the exact triangle \(F(A)\to F(B)\to
    F(C)\to F(\Sigma A)\) in~\(\Tri'\) splits.

  \end{itemize}
\end{lemma}

We will explain the notation during the proof.

\begin{proof}
  A morphism \(f\colon X\to Y\) in~\(\Tri'\) is called \emph{split epic}
  (\emph{split monic}) if there is \(g\colon Y\to X\) with \(f\circ
  g=\ID_Y\) (\(g\circ f=\ID_X\)).  An exact triangle \(X\overset{f}\to
  Y\overset{g}\to Z\overset{h}\to \Sigma X\) is said to \emph{split} if
  \(h=0\).  This immediately yields the characterisation of \(\ker
  F\)-exact triangles.  Any split triangle is isomorphic to a
  direct sum triangle, so that~\(f\) is split monic and~\(g\) is split
  epic (\cite{Neeman:Triangulated}*{Corollary 1.2.7}).  Conversely,
  either of these conditions implies that the triangle is split.

  Since the \(\ker F\)-exact triangles determine the \(\ker
  F\)-epimorphisms and \(\ker F\)-monomorphisms, the latter
  are detected by \(F(f)\) being split epic or split monic,
  respectively.  It is clear that split epimorphisms and split
  monomorphisms are epimorphisms and monomorphisms, respectively.  The
  converse holds in a triangulated category because if we embed a
  monomorphism or epimorphism in an exact triangle, then one of the maps
  is forced to vanish, so that the exact triangle splits.

  Finally, a morphism is invertible if and only if it is both split
  monic and split epic, and the zero map \(F(A)\to F(A)\) is invertible
  if and only if \(F(A)=0\).
\end{proof}

Alternatively, we may prove Lemma~\ref{lem:monic_kerF} using the Yoneda
embedding \(\Yoneda\colon \Tri'\to\Coh(\Tri')\).  The assertions about
phantom maps, equivalences, and contractibility boil down to the
observation that~\(\Yoneda\) is fully faithful.  The assertions about
monomorphisms and epimorphisms follow because a map \(f\colon A\to B\)
in~\(\Tri'\) becomes epic (monic) in \(\Coh(\Tri')\) if and only if it
is split epic (monic) in~\(\Tri'\).

There is a similar description for \(\bigcap \ker F_i\) for a set
\(\{F_i\}\) of exact functors.  This applies to the ideal~\(\VC_\Fam\)
for a family of (quantum) subgroups~\(\Fam\) in a locally compact
(quantum) group~\(G\) (Example~\ref{exa:VC}).  Replacing \(\VC_\Fam\)
by~\(\Fam\) in our notation to avoid clutter, we get:
\begin{itemize}
\item A morphism \(f\in\KK^G(A,B)\) is
  \begin{itemize}
  \item an \emph{\(\Fam\)\nb-phantom map} if and only if
    \(\Res_G^H(f)=0\) in~\(\KK^H\) for all \(H\in\Fam\);

  \item \emph{\(\Fam\)\nb-epic} if and only if \(\Res_G^H(f)\) is
    (split) epic in~\(\KK^H\) for all \(H\in\Fam\);

  \item \emph{\(\Fam\)\nb-monic} if and only if \(\Res_G^H(f)\) is
    (split) monic in~\(\KK^H\) for all \(H\in\Fam\);

  \item an \emph{\(\Fam\)\nb-equivalence} if and only if
    \(\Res_G^H(f)\) is a \(\KK^H\)\nb-equivalence for all
    \(H\in\Fam\).
  \end{itemize}

\item A \(G\)\nb-\(\Cst\)-algebra \(A\inOb\KK^G\) is
  \emph{\(\Fam\)\nb-contractible} if and only if \(\Res_G^H(A)\cong0\)
  in~\(\KK^H\) for all \(H\in\Fam\).

\item An exact triangle \(A\to B\to C\to \Sigma A\) in \(\KK^G\) is
  \emph{\(\Fam\)\nb-exact} if and only if
  \[
  \Res_G^H(A) \to \Res_G^H(B) \to \Res_G^H(C) \to \Sigma
  \Res_G^H(A)
  \]
  is a split exact triangle in~\(\KK^H\) for all \(H\in\Fam\).

\end{itemize}

You may write down a similar list for the ideal
\(\Ideal_\ltimes\subseteq\KK^G\) of Example~\ref{exa:cross}.

Lemma~\ref{lem:epic_kerF_exact} allows us to prove that the ideal
\(\IdealH\) in \(\Der(\Ab)\) cannot be the kernel of an exact functor:

\begin{proof}[Proof of~\textup{\ref{pro:counter_kernel_exact}}]
  We embed \(\Ab\to\Der(\Ab)\) as chain complexes concentrated in
  degree~\(0\).  The generator \(\tau\in\Ext(\Z/2,\Z/2)\) corresponds to
  the extension of Abelian groups \(\Z/2\into \Z/4\prto \Z/2\), where
  the first map is multiplication by~\(2\) and the second map is the
  natural projection.  We get an exact triangle
  \[
  \Z/2\to \Z/4 \to \Z/2 \xrightarrow{\tau} \Z/2[1]
  \]
  in \(\Der(\Ab)\).  This triangle is \(\IdealH\)\nb-exact because
  the map \(\Z/2\to \Z/4\) is injective as a group homomorphism and
  hence \(\IdealH\)\nb-monic in \(\Der(\Ab)\).

  Assume there were an exact functor \(F\colon \Der(\Ab)\to\Tri'\) with
  \(\ker F=\IdealH\).  Then \(F(\tau)=0\), so that~\(F\) maps our
  triangle to a split triangle and \(F(\Z/4)\cong F(\Z/2)\oplus
  F(\Z/2)\) by Lemma~\ref{lem:epic_kerF_exact}.  It follows that
  \(F(2\cdot \ID_{\Z/4}) = 2\cdot\ID_{F(\Z/4)}=0\) because \(2\cdot
  \ID_{F(\Z/2)}= F(2\cdot \ID_{\Z/2})=0\).  Hence \(2\cdot\ID_{\Z/4}\in
  \ker F=\IdealH\), which is false.  This contradiction shows that
  there is no exact functor~\(F\) with \(\ker F=\IdealH\).
\end{proof}

One of the most interesting questions about an ideal is whether all
\(\Ideal\)\nb-contractible objects vanish or, equivalently, whether all
\(\Ideal\)\nb-equivalences are invertible.  These two questions are
equivalent by Lemma~\ref{lem:equivalence_contractible}.  The answer is
negative for the ideal~\(\Ideal_\K\) in~\(\KK\) because the Universal
Coefficient Theorem does not hold for arbitrary separable
\(\Cst\)\nb-algebras.  Therefore, we also get counterexamples for the
ideal~\(\Ideal_{\cross,\K}\) in~\(\KK^G\) for a compact quantum group.
In contrast, if~\(G\) is a connected Lie group with torsion-free
fundamental group, then \(\Ideal_\ltimes\)\nb-equivalences in~\(\KK^G\)
are invertible (see~\cite{Meyer-Nest:BC_Coactions}).  If~\(G\) is an
amenable group, then \(\VC\)-equivalences in~\(\KK^G\) are invertible;
this follows from the proof of the Baum--Connes conjecture for these
groups by Nigel Higson and Gennadi Kasparov (see~\cite{Meyer-Nest:BC}).
These examples show that this question is subtle and may involve
difficult analysis.

\subsection{Exact chain complexes}
\label{sec:exact_chain_complexes}

The notion of \(\Ideal\)\nb-exactness, which we have only defined for
exact triangles so far, will now be extended to chain complexes.  Our
definition differs from Beligiannis' one
(\cites{Asadollahi-Salarian:Gorenstein, Beligiannis:Relative}), which we
recall first.

Let~\(\Tri\) be a triangulated category and let~\(\Ideal\) be a
homological ideal in~\(\Tri\).

\begin{definition}
  \label{def:decomposable_chain_complex}
  A chain complex
  \[
  C_\bullet \defeq (\dotsb\to
  C_{n+1} \xrightarrow{d_{n+1}}
  C_n \xrightarrow{d_n}
  C_{n-1} \xrightarrow{d_{n-1}}
  C_{n-2} \to \dotsb)
  \]
  in~\(\Tri\) is called \emph{\(\Ideal\)\nb-decomposable} if there is
  a sequence of \(\Ideal\)\nb-exact triangles
  \[
  K_{n+1} \xrightarrow{g_{n}}
  C_n \xrightarrow{f_n}
  K_{n} \xrightarrow{h_{n}}
  \Sigma K_{n+1}
  \]
  with \(d_n = g_{n-1}\circ f_n\colon C_n\to C_{n-1}\).
\end{definition}

Such complexes are called \(\Ideal\)\nb-exact in
\cites{Asadollahi-Salarian:Gorenstein, Beligiannis:Relative}.  This
definition is inspired by the following well-known fact: a chain complex
over an Abelian category is exact if and only if it splits into short
exact sequences of the form \(K_n\into C_n\prto K_{n-1}\) as in
Definition~\ref{def:decomposable_chain_complex}.

We prefer another definition of exactness because we have not found a
general explicit criterion for a chain complex to be
\(\Ideal\)\nb-decomposable.

\begin{definition}
  \label{def:exact_chain_complex}
  Let \(C_\bullet=(C_n,d_n)\) be a chain complex over~\(\Tri\).  For
  each \(n\in\N\), embed~\(d_n\) in an exact triangle
  \begin{equation}
    \label{eq:cone_over_differential}
    C_n \xrightarrow{d_n}
    C_{n-1} \xrightarrow{f_n}
    X_n \xrightarrow{g_n}
    \Sigma C_n.
  \end{equation}
  We call~\(C_\bullet\) \emph{\(\Ideal\)\nb-exact in degree~\(n\)} if
  the map \(X_{n} \xrightarrow{g_{n}}\Sigma C_{n} \xrightarrow{\Sigma
    f_{n+1}} \Sigma X_{n+1}\) belongs to \(\Ideal(X_n,\Sigma X_{n+1})\).
  This does not depend on auxiliary choices because the exact triangles
  in~\eqref{eq:cone_over_differential} are unique up to (non-canonical)
  isomorphism.

  We call~\(C_\bullet\) \emph{\(\Ideal\)\nb-exact} if it is
  \(\Ideal\)\nb-exact in degree~\(n\) for all \(n\in\Z\).
\end{definition}

This definition is designed to make the following lemma true:

\begin{lemma}
  \label{lem:exact_kerF}
  Let \(F\colon \Tri\to\Cat\) be a stable homological functor into a
  stable Abelian category~\(\Cat\) with \(\ker F=\Ideal\).  A chain
  complex~\(C_\bullet\) over~\(\Tri\) is \(\Ideal\)\nb-exact in
  degree~\(n\) if and only if
  \[
  F(C_{n+1}) \xrightarrow{F(d_{n+1})}
  F(C_n) \xrightarrow{F(d_n)}
  F(C_{n-1})
  \]
  is exact at \(F(C_n)\).
\end{lemma}

\begin{proof}
  The complex~\(C_\bullet\) is \(\Ideal\)-exact in degree~\(n\) if
  and only if the map
  \[
  \Sigma^{-1} F(X_{n}) \xrightarrow{\Sigma^{-1} F(g_{n})}
  F(C_{n}) \xrightarrow{F(f_{n+1})}
  F(X_{n+1})
  \]
  vanishes.  Equivalently, the range of \(\Sigma^{-1} F(g_{n})\) is
  contained in the kernel of \(F(f_{n+1})\).  The long exact sequences
  \begin{gather*}
    \dotsb \to
    \Sigma^{-1} F(X_{n}) \xrightarrow{\Sigma^{-1} F(g_{n})}
    F(C_{n}) \xrightarrow{F(d_{n})}
    F(C_{n-1}) \to \dotsb,\\
    \dotsb \to
    F(C_{n+1}) \xrightarrow{F(d_{n+1})}
    F(C_{n}) \xrightarrow{F(f_{n+1})}
    F(X_{n+1}) \to \dotsb
  \end{gather*}
  show that the range of \(\Sigma^{-1} F(g_{n})\) and the kernel of
  \(F(f_{n+1})\) are equal to the kernel of \(F(d_n)\) and the range of
  \(F(d_{n+1})\), respectively.  Hence~\(C_\bullet\) is
  \(\Ideal\)\nb-exact in degree~\(n\) if and only if \(\ker F(d_n)
  \subseteq \range F(d_{n+1})\).  Since \(d_n\circ d_{n+1}=0\), this is
  equivalent to \(\ker F(d_n) = \range F(d_{n+1})\).
\end{proof}

\begin{corollary}
  \label{cor:decomposable_exact}
  \(\Ideal\)\nb-decomposable chain complexes are
  \(\Ideal\)\nb-exact.
\end{corollary}

\begin{proof}
  Let \(F\colon \Tri\to\Cat\) be a stable homological functor with
  \(\ker F=\Ideal\).  If~\(C_\bullet\) is \(\Ideal\)\nb-decomposable,
  then \(F(C_\bullet)\) is obtained by splicing short exact sequences
  in~\(\Cat\).  This implies that \(F(C_\bullet)\) is exact, so
  that~\(C_\bullet\) is \(\Ideal\)\nb-exact by
  Lemma~\ref{lem:exact_kerF}.
\end{proof}

The converse implication in Corollary~\ref{cor:decomposable_exact} fails
in general (see Example~\ref{exa:exact_not_decomposable}).  But it holds
for chain complexes of projective objects; this follows from the proof
of Proposition~\ref{pro:resolutions}.

\begin{example}
  \label{exa:exact_chain}
  For the ideal~\(\Ideal_\K\) in~\(\KK\), Lemma~\ref{lem:exact_kerF}
  yields that a chain complex~\(C_\bullet\) over~\(\KK\) is
  \emph{\(\K\)\nb-exact} (in degree~\(n\)) if and only if the chain
  complex
  \[
  \dotsb\to \K_*(C_{n+1})\to \K_*(C_n)\to \K_*(C_{n-1})\to \dotsb
  \]
  of \(\Ztwo\)\nb-graded Abelian groups is exact (in degree~\(n\)).
  Similar remarks apply to the other ideals in~\S\ref{sec:examples} that
  are defined as kernels of stable homological functors.

  As a trivial example, we consider the largest possible ideal
  \(\Ideal=\Tri\).  This ideal is defined by the zero functor.
  Lemma~\ref{lem:exact_kerF} or the definition yield that \emph{all}
  chain complexes are \(\Tri\)\nb-exact.  In contrast, it seems hard
  to characterise the \(\Ideal\)\nb-decomposable chain complexes,
  already for \(\Ideal=\Tri\).
\end{example}

\begin{lemma}
  \label{lem:exact_complex_length_3}
  A chain complex of length~\(3\)
  \[
  \dotsb \to 0 \to
  A \xrightarrow{f}
  B \xrightarrow{g}
  C \to
  0 \to \dotsb
  \]
  is \(\Ideal\)\nb-exact if and only if there are an
  \(\Ideal\)\nb-exact exact triangle \(A'\xrightarrow{f'}
  B'\xrightarrow{g'} C'\to \Sigma A'\) and a commuting diagram
  \begin{equation}
    \label{eq:exact_complex_length_3}
    \begin{gathered}
      \xymatrix{
        A' \ar[r]^{f'} \ar[d]^{\alpha}_{\sim} &
        B' \ar[r]^{g'} \ar[d]^{\beta}_{\sim} &
        C' \ar[d]^{\gamma}_{\sim} \\
        A \ar[r]^{f} &
        B \ar[r]^{g} &
        C
      }
    \end{gathered}
  \end{equation}
  where the vertical maps \(\alpha,\beta,\gamma\) are
  \(\Ideal\)\nb-equivalences.  Furthermore, we can achieve that
  \(\alpha\) and~\(\beta\) are identity maps.
\end{lemma}

\begin{proof}
  Let~\(F\) be a stable homological functor with \(\Ideal=\ker F\).

  Suppose first that we are in the situation
  of~\eqref{eq:exact_complex_length_3}.  Lemma~\ref{lem:epic_kerF}
  yields that \(F(\alpha)\), \(F(\beta)\), and \(F(\gamma)\) are
  invertible and that \(0\to F(A')\to F(B')\to F(C')\to 0\) is a short
  exact sequence.  Hence so is \(0\to F(A)\to F(B)\to F(C)\to 0\).  Now
  Lemma~\ref{lem:exact_kerF} yields that our given chain complex is
  \(\Ideal\)\nb-exact.

  Conversely, suppose that we have an \(\Ideal\)\nb-exact chain
  complex.  By Lemma~\ref{lem:exact_kerF}, this means that \(0\to
  F(A)\to F(B)\to F(C)\to 0\) is a short exact sequence.  Hence
  \(f\colon A\to B\) is \(\Ideal\)\nb-monic.  Embed~\(f\) in an exact
  triangle \(A\to B\to C'\to \Sigma A\).  Since~\(f\) is
  \(\Ideal\)\nb-monic, this triangle is \(\Ideal\)\nb-exact.  Let
  \(\alpha=\ID_A\) and \(\beta=\ID_B\).  Since the functor
  \(\Tri(\blank,C)\) is cohomological and \(g\circ f=0\), we can find a
  map \(\gamma\colon C'\to C\) making~\eqref{eq:exact_complex_length_3}
  commute.  The functor~\(F\) maps the rows
  of~\eqref{eq:exact_complex_length_3} to short exact sequences by
  Lemmas \ref{lem:exact_kerF} and~\ref{lem:epic_kerF}.  Now the Five
  Lemma yields that \(F(\gamma)\) is invertible, so that~\(\gamma\) is
  an \(\Ideal\)\nb-equivalence.
\end{proof}

\begin{remark}
  \label{rem:length_3_exact_decomposable}
  Lemma~\ref{lem:exact_complex_length_3} implies that
  \(\Ideal\)\nb-exact chain complexes of length~\(3\) are
  \(\Ideal\)\nb-decomposable.  We do not expect this for chain
  complexes of length~\(4\).  But we have not searched for a
  counterexample.
\end{remark}

Which chain complexes over~\(\Tri\) are \(\Ideal\)\nb-exact for
\(\Ideal=0\) and hence for any homological ideal?  The next definition
provides the answer.

\begin{definition}
  \label{def:homologically_exact_chain_complex}
  A chain complex~\(C_\bullet\) over a triangulated category is called
  \emph{homologically exact} if \(F(C_\bullet)\) is exact for any
  homological functor \(F\colon \Tri\to\Cat\).
\end{definition}

\begin{example}
  \label{exa:exact_triangle_homologically_exact}
  If \(A\to B\to C\to\Sigma A\) is an exact triangle, then the chain
  complex
  \[
  \dotsb \to \Sigma^{-1} A \to \Sigma^{-1} B \to \Sigma^{-1} C \to A
  \to B \to C \to \Sigma A\to \Sigma B\to \Sigma C \to \dotsb
  \]
  is homologically exact by the definition of a homological functor.
\end{example}

\begin{lemma}
  \label{lem:homologically_exact_chain_complexes}
  Let \(F\colon \Tri\to\Tri'\) be an exact functor between two
  triangulated categories.  Let~\(C_\bullet\) be a chain complex
  over~\(\Tri\).  The following are equivalent:
  \begin{enumerate}[leftmargin=*,label=\textup{(\arabic*)}]
  \item \(C_\bullet\) is \(\ker F\)-exact in degree~\(n\);

  \item \(F(C_\bullet)\) is \(\Ideal\)\nb-exact in degree~\(n\) with
    respect to \(\Ideal=0\);

  \item the chain complex \(\Yoneda\circ F(C_\bullet)\) in
    \(\Coh(\Tri')\) is exact in degree~\(n\);

  \item \(F(C_\bullet)\) is homologically exact in degree~\(n\);

  \item the chain complexes of Abelian groups
    \(\Tri'\bigl(A,F(C_\bullet)\bigr)\) are exact in degree~\(n\) for
    all \(A\inOb\Tri'\).

  \end{enumerate}
\end{lemma}

\begin{proof}
  By Freyd's Theorem~\ref{the:Freyd}, \(\Yoneda\circ F\colon
  \Tri\to\Coh(\Tri')\) is a stable homological functor with \(\ker
  F=\ker (\Yoneda\circ F)\).  Hence Lemma~\ref{lem:exact_kerF} yields
  (1)\(\iff\)(3).  Similarly, we have (2)\(\iff\)(3) because
  \(\Yoneda\colon \Tri'\to\Coh(\Tri')\) is a stable homological functor
  with \(\ker\Yoneda=0\).  Freyd's Theorem~\ref{the:Freyd} also asserts
  that any homological functor \(F\colon \Tri'\to\Cat'\) factors as
  \(\bar{F}\circ\Yoneda\) for an exact functor~\(\bar{F}\).  Hence
  (3)\(\Longrightarrow\)(4).
  Proposition~\ref{pro:triangulated_rep_cohomological} yields
  (4)\(\Longrightarrow\)(5).  Finally, (5)\(\iff\)(3) because kernels
  and cokernels in \(\Coh(\Tri')\) are computed pointwise on objects
  of~\(\Tri'\).
\end{proof}

\begin{remark}
  \label{rem:intersection}
  More generally, consider a set of exact functors \(F_i\colon
  \Tri\to\Tri'_i\).  As in the proof of the equivalence (1)\(\iff\)(2)
  in Lemma~\ref{lem:homologically_exact_chain_complexes}, we see that a
  chain complex~\(C_\bullet\) is \(\bigcap \ker F_i\)-exact (in
  degree~\(n\)) if and only if the chain complexes \(F_i(C_\bullet)\)
  are exact (in degree~\(n\)) for all~\(i\).
\end{remark}

As a consequence, a chain complex~\(C_\bullet\) over~\(\KK^G\) for a
locally compact quantum group~\(G\) is \(\Fam\)\nb-exact if and only
if \(\Res_G^H(C_\bullet)\) is homologically exact for all \(H\in\Fam\).
A chain complex~\(C_\bullet\) over~\(\KK^G\) for a compact quantum
group~\(G\) is \(\Ideal_\ltimes\)-exact if and only if the chain
complex \(G\cross C_\bullet\) over~\(\KK\) is homologically exact.

\begin{example}
  \label{exa:exact_not_decomposable}
  We exhibit an \(\Ideal\)\nb-exact chain complex that is not
  \(\Ideal\)\nb-decomposable for the ideal \(\Ideal=0\).  By
  Lemma~\ref{lem:epic_kerF_exact}, any \(0\)\nb-exact triangle is
  split.  Therefore, a chain complex is \(0\)\nb-decomposable if and
  only if it is a direct sum of chain complexes of the form \(0 \to
  K_n\xrightarrow{\ID} K_n\to 0\).  Hence any decomposable chain complex
  is contractible and therefore mapped by any homological functor to a
  contractible chain complex.  By the way, if idempotents in~\(\Tri\)
  split then a chain complex is \(0\)\nb-decomposable if and only if
  it is contractible.

  As we have remarked in
  Example~\ref{exa:exact_triangle_homologically_exact}, the chain
  complex
  \[
  \dotsb \to \Sigma^{-1} C \to A\to B\to C \to \Sigma A\to \Sigma B\to \Sigma
  C \to \Sigma^2 A\to \dotsb
  \]
  is homologically exact for any exact triangle \(A\to B\to C\to\Sigma
  A\).  But such chain complexes need not be contractible.  A
  counterexample is the exact triangle \(\Ztwo\to \Z/4\to\Ztwo\to
  \Sigma\Ztwo\) in \(\Der(\Ab)\), which we have already used in the
  proof of Proposition~\ref{pro:counter_kernel_exact}.  The resulting
  chain complex over \(\Der(\Ab)\) cannot be contractible
  because~\(\Hgy\) maps it to a non-contractible chain complex.
\end{example}

\subsubsection{More homological algebra with chain complexes}
\label{sec:homological_algebra_chain}

Using our notion of exactness for chain complexes, we can do homological
algebra in the homotopy category \(\Ho(\Tri)\).  We briefly sketch some
results in this direction, assuming some familiarity with more advanced
notions from homological algebra.  We will not use this later.

The \(\Ideal\)\nb-exact chain complexes form a thick subcategory of
\(\Ho(\Tri)\) because of Lemma~\ref{lem:exact_kerF}.  We let
\(\Der\defeq \Der(\Tri,\Ideal)\) be the localisation of \(\Ho(\Tri)\) at
this subcategory and call it the \emph{derived category of~\(\Tri\) with
  respect to~\(\Ideal\)}.

We let \(\Der^{\ge n}\) and \(\Der^{\le n}\) be the full subcategories
of \(\Der\) consisting of chain complexes that are \(\Ideal\)\nb-exact
in degrees less than~\(n\) and greater than~\(n\), respectively.

\begin{theorem}
  \label{the:t-structure}
  The pair of subcategories \(\Der^{\ge0}\), \(\Der^{\le0}\) forms a
  \emph{truncation structure} \textup{(}t-structure\textup{)} on
  \(\Der\) in the sense of~\cite{Beilinson-Bernstein-Deligne}.
\end{theorem}

\begin{proof}
  The main issue here is the truncation of chain complexes.
  Let~\(C_\bullet\) be a chain complex over~\(\Tri\).  We embed the
  map~\(d_0\) in an exact triangle \(C_0 \to C_{-1} \to X\to \Sigma
  C_0\) and let \(C_\bullet^{\ge0}\) be the chain complex
  \[
  \dotsb \to C_2 \to C_1 \to C_0 \to C_{-1}
  \to X \to \Sigma C_0 \to \Sigma C_{-1} \to \Sigma X \to \Sigma^2 C_0
  \to \dotsb.
  \]
  This chain complex is \(\Ideal\)\nb-exact -- even homologically exact
  -- in negative degrees, that is, \(C_\bullet^{\ge0}\in\Der^{\ge0}\).
  The triangulated category structure allows us to construct a chain map
  \(C_\bullet^{\ge0}\to C_\bullet\) that is an isomorphism on~\(C_n\)
  for \(n\ge-1\).  Hence its mapping cone \(C_\bullet^{\le-1}\) is
  \(\Ideal\)\nb-exact -- even contractible -- in
  degrees~\({}\ge0\), that is, \(C_\bullet^{\le-1}\inOb\Der^{\le-1}\).
  By construction, we have an exact triangle
  \[
  C_\bullet^{\ge0}\to C_\bullet\to C_\bullet^{\le-1}\to
  \Sigma C_\bullet^{\ge0}
  \]
  in \(\Der\).

  We also have to check that there is no non-zero morphism
  \(C_\bullet\to D_\bullet\) in \(\Der\) if
  \(C_\bullet\inOb\Der^{\ge0}\) and \(D_\bullet\inOb\Der^{\le-1}\).
  Recall that morphisms in \(\Der\) are represented by diagrams
  \(C_\bullet \overset{\sim}\leftarrow \tilde{C}_\bullet \to D_\bullet\)
  in \(\Ho(\Tri)\), where the first map is an
  \(\Ideal\)\nb-equivalence.  Hence \(\tilde{C}_\bullet\inOb
  \Der^{\ge0}\) as well.  We claim that any chain map \(f\colon
  \tilde{C}_\bullet^{\ge0} \to D_\bullet^{\le-1}\) is homotopic
  to~\(0\).  Since the maps \(\tilde{C}_\bullet^{\ge0} \to C_\bullet\)
  and \(D_\bullet \to D_\bullet^{\le-1}\) are
  \(\Ideal\)\nb-equivalences, any morphism \(C_\bullet\to D_\bullet\)
  vanishes in \(\Der\).

  It remains to prove the claim.  In a first step, we use that
  \(D_\bullet^{\le-1}\) is contractible in degrees~\({}\ge0\) to
  replace~\(f\) by a homotopic chain map supported in degrees~\(<0\).
  In a second step, we use that \(\tilde{C}_\bullet^{\ge0}\) is
  homologically exact in the relevant degrees to recursively construct a
  chain homotopy between~\(f\) and~\(0\).
\end{proof}

Any truncation structure gives rise to an Abelian category, its
\emph{core}.  The core of~\(\Der\) is the full subcategory~\(\Cat\) of
all chain complexes that are \(\Ideal\)\nb-exact except in degree~\(0\).
This is a stable Abelian category, and the standard embedding
\(\Tri\to\Ho(\Tri)\) yields a stable homological functor \(F\colon
\Tri\to\Cat\) with \(\ker F=\Ideal\).

This functor is characterised uniquely by the following universal
property: any (stable) homological functor \(H\colon \Tri\to\Cat'\) with
\(\Ideal\subseteq \ker H\) factors uniquely as \(H=\bar{H}\circ F\) for
an exact functor \(\bar{H}\colon \Cat\to\Cat'\).  We
construct~\(\bar{H}\) in three steps.

First, we lift~\(H\) to an exact functor \(\Ho(H)\colon
\Ho(\Tri,\Ideal)\to\Ho(\Cat')\).  Secondly, \(\Ho(H)\) descends to a
functor \(\Der(H)\colon \Der(\Tri,\Ideal)\to\Der(\Cat')\).  Finally,
\(\Der(H)\) restricts to a functor \(\bar{H}\colon \Cat\to\Cat'\)
between the cores.  Since \(\Ideal\subseteq \ker H\), an
\(\Ideal\)\nb-exact chain complex is also \(\ker H\)\nb-exact.
Hence \(\Ho(H)\) preserves exactness of chain complexes by
Lemma~\ref{lem:exact_kerF}.  This allows us to construct \(\Der(H)\) and
shows that \(\Der(H)\) is compatible with truncation structures.  This
allows us to restrict it to an exact functor between the cores.
Finally, we use that the core of the standard truncation structure on
\(\Der(\Cat)\) is~\(\Cat\).  It is easy to see that we have
\(\bar{H}\circ F=H\).

Especially, we get an exact functor \(\Der(F)\colon
\Der(\Tri,\Ideal)\to\Der(\Cat)\), which restricts to the identity
functor \(\ID_\Cat\) between the cores.  Hence \(\Der(F)\) is fully
faithful on the thick subcategory generated by
\(\Cat\subseteq\Der(\Tri,\Ideal)\).  It seems plausible that \(\Der(F)\)
should be an equivalence of categories under some mild conditions
on~\(\Ideal\) and~\(\Tri\).

We will continue our study of the functor \(F\colon \Tri\to\Cat\)
in~\S\ref{sec:universal_functor}.  The universal property determines it
uniquely.  Beligiannis (\cite{Beligiannis:Relative}) has another,
simpler construction.

\subsection{Projective objects}
\label{sec:projective}

Let~\(\Ideal\) be a homological ideal in a triangulated
category~\(\Tri\).

\begin{definition}
  \label{def:ideal_exact_functor}\label{def:ideal_projective}\label{def:ideal_injective}
  A homological functor \(F\colon \Tri\to\Cat\) is called
  \emph{\(\Ideal\)\nb-exact} if \(F(f)=0\) for all
  \(\Ideal\)\nb-phantom maps~\(f\) or, equivalently, \(\Ideal\subseteq
  \ker F\).  An object \(A\inOb\Tri\) is called
  \emph{\(\Ideal\)\nb-projective} if the functor
  \(\Tri(A,\blank)\colon \Tri\to\Ab\) is \(\Ideal\)\nb-exact.  Dually,
  an object \(B\inOb\Tri\) is called \emph{\(\Ideal\)\nb-injective} if
  the functor \(\Tri(\blank,B)\colon \Tri\to\Ab^\op\) is
  \(\Ideal\)\nb-exact.

  We write~\(\Proj_\Ideal\) for the class of \(\Ideal\)\nb-projective
  objects in~\(\Tri\).
\end{definition}

The notions of projective and injective object are dual to each other:
if we pass to the opposite category~\(\Tri^\op\) with the canonical
triangulated category structure and use the same ideal~\(\Ideal^\op\),
then this exchanges the roles of projective and injective objects.
Therefore, it suffices to discuss one of these two notions in the
following.  We will only treat projective objects because all the ideals
in~\S\ref{sec:examples} have enough projective objects, but most of them
do not have enough injective objects.

Notice that the functor~\(F\) is \(\Ideal\)\nb-exact if and only if
the associated stable functor \(F_*\colon \Tri\to\Cat^\Z\) is
\(\Ideal\)\nb-exact because~\(\Ideal\) is stable.

Since we require~\(F\) to be homological, the long exact homology
sequence and Lemma~\ref{lem:exact_kerF} yield that the following
conditions are all equivalent to~\(F\) being \(\Ideal\)\nb-exact:
\begin{itemize}
\item \(F\) maps \(\Ideal\)\nb-epimorphisms to epimorphisms
  in~\(\Cat\);

\item \(F\) maps \(\Ideal\)\nb-monomorphisms to monomorphisms
  in~\(\Cat\);

\item \(0\to F(A)\to F(B)\to F(C)\to0\) is a short exact sequence
  in~\(\Cat\) for any \(\Ideal\)\nb-exact triangle \(A\to B\to
  C\to\Sigma A\);

\item \(F\) maps \(\Ideal\)\nb-exact chain complexes to exact chain
  complexes in~\(\Cat\).

\end{itemize}
This specialises to equivalent definitions of \(\Ideal\)\nb-projective
objects.

\begin{lemma}
  \label{lem:ideal_projective_injective}
  An object \(A\inOb\Tri\) is \(\Ideal\)\nb-projective if and only if
  \(\Ideal(A,B)=0\) for all \(B\inOb\Tri\).
\end{lemma}

\begin{proof}
  If \(f\in\Ideal(A,B)\), then \(f=f_*(\ID_A)\).  This has to vanish
  if~\(A\) is \(\Ideal\)\nb-projective.  Suppose, conversely, that
  \(\Ideal(A,B)=0\) for all \(B\inOb\Tri\).  If \(f\in\Ideal(B,B')\),
  then \(\Tri(A,f)\) maps \(\Tri(A,B)\) to \(\Ideal(A,B')=0\), so that
  \(\Tri(A,f)=0\).  Hence~\(A\) is \(\Ideal\)\nb-projective.
\end{proof}

An \(\Ideal\)\nb-exact functor also has the following properties
(which are strictly weaker than being \(\Ideal\)\nb-exact):
\begin{itemize}
\item \(F\) maps \(\Ideal\)\nb-equivalences to isomorphisms
  in~\(\Cat\);

\item \(F\) maps \(\Ideal\)\nb-contractible objects to~\(0\)
  in~\(\Cat\).

\end{itemize}
Again we may specialise this to \(\Ideal\)\nb-projective objects.

\begin{lemma}
  \label{lem:projective_hereditary}\label{lem:injective_hereditary}
  The class of \(\Ideal\)\nb-exact homological functors \(\Tri\to\Ab\)
  or \(\Tri\to\Ab^\op\) is closed under composition with
  \(\Sigma^{\pm1}\colon \Tri\to\Tri\), retracts, direct sums, and direct
  products.  The class~\(\Proj_\Ideal\) of \(\Ideal\)\nb-projective
  objects is closed under \textup{(}de\textup{)}sus\-pen\-sions,
  retracts, and possibly infinite direct sums \textup{(}as far as they
  exist in~\(\Tri\)\textup{)}.
\end{lemma}

\begin{proof}
  The first assertion follows because direct sums and products of
  Abelian groups are exact; the second one is a special case.
\end{proof}

\begin{notation}
  \label{note:sums_closure}
  Let~\(\Proj\) be a set of objects of~\(\Tri\).  We let
  \(\sumclosure{\Proj}\) be the smallest class of objects of~\(\Tri\)
  that contains~\(\Proj\) and is closed under retracts and direct sums
  (as far as they exist in~\(\Tri\)).
\end{notation}

By Lemma~\ref{lem:projective_hereditary}, \(\sumclosure{\Proj}\)
consists of \(\Ideal\)\nb-projective objects if~\(\Proj\) does.  We
say that~\(\Proj\) \emph{generates all \(\Ideal\)\nb-projective
  objects} if \(\sumclosure{\Proj} = \Proj_\Ideal\).  In examples, it is
usually easier to describe a class of generators in this sense.

\subsection{Projective resolutions}
\label{sec:proj_resolutions}

\begin{definition}
  \label{def:one-step_projective_resolution}
  Let \(\Ideal\subseteq\Tri\) be a homological ideal in a triangulated
  category and let \(A\inOb\Tri\).  A \emph{one-step
    \(\Ideal\)\nb-projective resolution} is an
  \(\Ideal\)\nb-epimorphism \(\pi\colon P\to A\) with
  \(P\inOb\Proj_\Ideal\).  An \emph{\(\Ideal\)\nb-projective
    resolution of~\(A\)} is an \(\Ideal\)\nb-exact chain complex
  \[
  \dotsb \xrightarrow{\delta_{n+1}}
  P_n \xrightarrow{\delta_n}
  P_{n-1} \xrightarrow{\delta_{n-1}}
  \dotsb \xrightarrow{\delta_1}
  P_0 \xrightarrow{\delta_0} A
  \]
  with \(P_n\inOb\Proj_\Ideal\) for all \(n\in\N\).

  We say that~\(\Ideal\) has \emph{enough projective objects} if each
  \(A\inOb\Tri\) has a one-step \(\Ideal\)\nb-projective resolution.
\end{definition}

The following proposition contains the basic properties of projective
resolutions, which are familiar from the similar situation for Abelian
categories.

\begin{proposition}
  \label{pro:resolutions}
  If~\(\Ideal\) has enough projective objects, then any object
  of~\(\Tri\) has an \(\Ideal\)\nb-projective resolution \textup{(}and
  vice versa\textup{)}.

  Let \(P_\bullet\to A\) and \(P_\bullet'\to A'\) be
  \(\Ideal\)\nb-projective resolutions.  Then any map \(A\to A'\) may
  be lifted to a chain map \(P_\bullet\to P_\bullet'\), and this lifting
  is unique up to chain homotopy.  Two \(\Ideal\)\nb-projective
  resolutions of the same object are chain homotopy equivalent.  As a
  result, the construction of projective resolutions provides a functor
  \[
  P\colon \Tri\to \Ho(\Tri).
  \]

  Let \(A\overset{f}\to B\overset{g}\to C\overset{h}\to \Sigma A\) be an
  \(\Ideal\)\nb-exact triangle.  Then there exists a canonical map
  \(\eta\colon P(C)\to P(A)[1]\) in \(\Ho(\Tri)\) such that the triangle
  \[
  P(A) \xrightarrow{P(f)}
  P(B) \xrightarrow{P(g)}
  P(C) \xrightarrow{\eta}
  P(A)[1]
  \]
  in \(\Ho(\Tri)\) is exact; here~\([1]\) denotes the translation
  functor in \(\Ho(\Tri)\), which has nothing to do with the suspension
  in~\(\Tri\).
\end{proposition}

\begin{proof}
  Let \(A\inOb\Tri\).  By assumption, there is a one-step
  \(\Ideal\)-projective resolution \(\delta_0\colon P_0\to A\),
  which we embed in an exact triangle \(A_1 \to P_0\to A\to \Sigma
  A_1\).  Since~\(\delta_0\) is \(\Ideal\)\nb-epic, this triangle is
  \(\Ideal\)\nb-exact.  By induction, we construct a sequence of such
  \(\Ideal\)\nb-exact triangles \(A_{n+1} \to P_n \to A_n \to \Sigma
  A_{n+1}\) for \(n\in\N\) with \(P_n\inOb\Proj\) and \(A_0=A\).  By
  composition, we obtain maps \(\delta_n\colon P_n\to P_{n-1}\) for
  \(n\ge1\), which satisfy \(\delta_n\circ\delta_{n+1}=0\) for all
  \(n\ge0\).  The resulting chain complex
  \[
  \dotsb \to
  P_n \xrightarrow{\delta_n}
  P_{n-1} \xrightarrow{\delta_{n-1}}
  P_{n-2} \to
  \dotsb \to
  P_1 \xrightarrow{\delta_1}
  P_0 \xrightarrow{\delta_0}
  A \to 0
  \]
  is \(\Ideal\)\nb-decomposable by construction and therefore
  \(\Ideal\)\nb-exact by Corollary~\ref{cor:decomposable_exact}.

  The remaining assertions are proved exactly as their classical
  counterparts in homological algebra.  We briefly sketch the arguments.
  Let \(P_\bullet \to A\) and \(P_\bullet'\to A'\) be
  \(\Ideal\)\nb-projective resolutions and let \(f\in \Tri(A,A')\).
  We construct \(f_n\in\Tri(P_n,P_n')\) by induction on~\(n\) such that
  the diagrams
  \[
  \begin{gathered}
    \xymatrix{
      P_0 \ar[r]^{\delta_0} \ar@{.>}[d]_{f_0} \ar[dr] &
      A \ar[d]^{f} \\
      P'_0 \ar[r]_{\delta'_0} & A' \\
    }
  \end{gathered}
  \qquad\text{and}\qquad
  \begin{gathered}
    \xymatrix{
      P_n \ar[r]^{\delta_n} \ar@{.>}[d]_{f_n} \ar[dr] &
      P_{n-1} \ar[d]^{f_{n-1}} \\
      P'_n \ar[r]_{\delta'_n} & P'_{n-1} \\
    }
  \end{gathered}
  \quad\text{for \(n\ge1\)}
  \]
  commute.  We must check that this is possible.  Since the chain
  complex \(P_\bullet'\to A\) is \(\Ideal\)\nb-exact and~\(P_n\) is
  \(\Ideal\)\nb-projective for all \(n\ge0\), the chain complexes
  \[
  \dotsb \to
  \Tri(P_n,P_m') \xrightarrow{(\delta_m')_*}
  \Tri(P_n,P_{m-1}') \to
  \dotsb \to
  \Tri(P_n,P_0') \xrightarrow{(\delta_0')_*}
  \Tri(P_n,A) \to 0
  \]
  are exact for all \(n\in\N\).  This allows us to find maps~\(f_n\) as
  above.  By construction, these maps form a chain map lifting \(f\colon
  A\to A'\).  Its uniqueness up to chain homotopy is proved similarly.
  If we apply this unique lifting result to two \(\Ideal\)\nb-projective
  resolutions of the same object, we get the uniqueness of
  \(\Ideal\)\nb-projective resolutions up to chain homotopy equivalence.
  Hence we get a well-defined functor \(P\colon \Tri\to\Ho(\Tri)\).

  Now consider an \(\Ideal\)\nb-exact triangle \(A\to B\to C\to \Sigma
  A\) as in the third paragraph of the lemma.  Let~\(X_\bullet\) be the
  mapping cone of some chain map \(P(A)\to P(B)\) in the homotopy class
  \(P(f)\).  This chain complex is supported in degrees~\({}\ge0\) and
  has \(\Ideal\)\nb-projective entries because \(X_n = P(A)_{n-1} \oplus
  P(B)_n\).  The map \(X_0=0\oplus P(B)_0\to B\to C\) yields a chain map
  \(X_\bullet\to C\), that is, the composite map \(X_1\to X_0\to C\)
  vanishes.  By construction, this chain map lifts the given map \(B\to
  C\) and we have an exact triangle \(P(A)\to P(B) \to X_\bullet\to
  P(A)[1]\) in \(\Ho(\Tri)\).  It remains to observe that \(X_\bullet\to
  C\) is \(\Ideal\)\nb-exact.  Then~\(X_\bullet\) is an
  \(\Ideal\)\nb-projective resolution of~\(C\).  Since such resolutions
  are unique up to chain homotopy equivalence, we get a canonical
  isomorphism \(X_\bullet\cong P(C)\) in \(\Ho(\Tri)\) and hence the
  assertion in the third paragraph.

  Let~\(F\) be a stable homological functor with \(\Ideal=\ker F\).  We
  have to check that \(F(X_\bullet)\to F(C)\) is a resolution.  This
  reduces to a well-known diagram chase in Abelian categories, using
  that \(F\bigl(P(A)\bigr)\to F(A)\) and \(F\bigl(P(B)\bigr)\to F(B)\)
  are resolutions and that \(F(A)\into F(B)\prto F(C)\) is exact.
\end{proof}

\subsection{Derived functors}
\label{sec:derived_functors}

We only define derived functors if there are enough projective objects
because this case is rather easy and suffices for our applications.  The
general case can be reduced to the familiar case of Abelian categories
using the results of~\S\ref{sec:homological_algebra_chain}.

\begin{definition}
  \label{def:derived}
  Let~\(\Ideal\) be a homological ideal in a triangulated
  category~\(\Tri\) with enough projective objects.  Let \(F\colon
  \Tri\to\Cat\) be an additive functor with values in an Abelian
  category~\(\Cat\).  It induces a functor \(\Ho(F)\colon \Ho(\Tri)\to
  \Ho(\Cat)\), applying~\(F\) pointwise to chain complexes.  Let
  \(P\colon \Tri\to\Ho(\Tri)\) be the projective resolution functor
  constructed in Proposition~\ref{pro:resolutions}.  Let \(H_n\colon
  \Ho(\Cat)\to\Cat\) be the \(n\)th homology functor for some
  \(n\in\N\).  The composite functor
  \[
  \Left_n F\colon \Tri \xrightarrow{P}
  \Ho(\Tri) \xrightarrow{\Ho(F)}
  \Ho(\Cat) \xrightarrow{H_n}
  \Cat
  \]
  is called the \(n\)th \emph{left derived functor} of~\(F\).  If
  \(F\colon \Tri^\op\to\Cat\) is an additive functor, then the
  corresponding functor \(H^n\circ \Ho(F)\circ P\colon \Tri^\op\to\Cat\)
  is denoted by \(\Right^n F\) and called the \(n\)th \emph{right
    derived functor} of~\(F\).
\end{definition}

More concretely, let \(A\inOb\Tri\) and let
\((P_\bullet,\delta_\bullet)\) be an \(\Ideal\)\nb-projective resolution
of~\(A\).  If~\(F\) is covariant, then \(\Left_n F(A)\) is the homology
at~\(F(P_n)\) of the chain complex
\[
\dotsb \to
F(P_{n+1}) \xrightarrow{F(\delta_{n+1})}
F(P_n) \xrightarrow{F(\delta_n)}
F(P_{n-1}) \to
\dotsb \to
F(P_0) \to 0.
\]
If~\(F\) is contravariant, then \(\Right^n F(A)\) is the cohomology at
\(F(P_n)\) of the cochain complex
\[
\dotsb \leftarrow
F(P_{n+1}) \xleftarrow{F(\delta_{n+1})}
F(P_n) \xleftarrow{F(\delta_n)}
F(P_{n-1}) \leftarrow
\dotsb \leftarrow
F(P_0)\leftarrow 0.
\]

\begin{lemma}
  \label{lem:derived_functors_les}
  Let \(A\to B\to C\to\Sigma A\) be an \(\Ideal\)\nb-exact triangle.
  If \(F\colon \Tri\to\Cat\) is a covariant additive functor, then there
  is a long exact sequence
  \begin{multline*}
    \dotsb
    \to \Left_n F(A)
    \to \Left_n F(B)
    \to \Left_n F(C)
    \to \Left_{n-1} F(A)\\
    \to \dotsb
    \to \Left_1F(C)
    \to \Left_0F(A)
    \to \Left_0F(B)
    \to \Left_0F(C)
    \to 0.
  \end{multline*}
  If~\(F\) is contravariant instead, then there is a long exact sequence
  \begin{multline*}
    \dotsb \leftarrow \Right^n F(A)
    \leftarrow \Right^n F(B)
    \leftarrow \Right^n F(C)
    \leftarrow \Right^{n-1} F(A)\\
    \leftarrow \dotsb
    \leftarrow \Right^1F(C)
    \leftarrow \Right^0F(A)
    \leftarrow \Right^0F(B)
    \leftarrow \Right^0F(C)
    \leftarrow 0.
  \end{multline*}
\end{lemma}

\begin{proof}
  This follows from the third assertion of
  Proposition~\ref{pro:resolutions} together with the well-known long
  exact homology sequence for exact triangles in \(\Ho(\Cat)\).
\end{proof}

\begin{lemma}
  \label{lem:characterise_exact_functors}
  Let \(F\colon \Tri\to\Cat\) be a homological functor.  The following
  assertions are equivalent:
  \begin{enumerate}[leftmargin=*,label=\textup{(\arabic*)}]
  \item \(F\) is \(\Ideal\)\nb-exact;

  \item \(\Left_0 F(A)\cong F(A)\) and \(\Left_p F(A)=0\) for all \(p>0\),
    \(A\inOb\Tri\);

  \item \(\Left_0 F(A)\cong F(A)\) for all \(A\inOb\Tri\).

  \end{enumerate}
  The analogous assertions for contravariant functors are equivalent as
  well.
\end{lemma}

\begin{proof}
  If~\(F\) is \(\Ideal\)\nb-exact, then~\(F\) maps
  \(\Ideal\)\nb-exact chain complexes in~\(\Tri\) to exact chain
  complexes in~\(\Cat\).  This applies to \(\Ideal\)\nb-projective
  resolutions, so that (1)\(\Longrightarrow\)(2)\(\Longrightarrow\)(3).
  It follows from (3) and Lemma~\ref{lem:derived_functors_les}
  that~\(F\) maps \(\Ideal\)\nb-epimorphisms to epimorphisms.  Since
  this characterises \(\Ideal\)\nb-exact functors, we get
  (3)\(\Longrightarrow\)(1).
\end{proof}

It can happen that \(\Left_p F=0\) for all \(p>0\) although~\(F\) is not
\(\Ideal\)\nb-exact.

We have a natural transformation \(\Left_0 F(A)\to F(A)\) (or \(F(A)\to
\Right^0 F(A)\)), which is induced by the augmentation map \(P_\bullet
\to A\) for an \(\Ideal\)\nb-projective resolution.
Lemma~\ref{lem:characterise_exact_functors} shows that these maps are
usually not bijective, although this happens frequently for derived
functors on Abelian categories.

\begin{definition}
  \label{def:Ext}
  We let \(\Ext^n_{\Tri,\Ideal}(A,B)\) be the \(n\)th right derived
  functor with respect to~\(\Ideal\) of the contravariant functor
  \(A\mapsto \Tri(A,B)\).
\end{definition}

We have natural maps \(\Tri(A,B)\to\Ext^0_{\Tri,\Ideal}(A,B)\), which
usually are not invertible.  Lemma~\ref{lem:derived_functors_les} yields
long exact sequences
\begin{multline*}
  \dotsb
  \leftarrow \Ext^n_{\Tri,\Ideal}(A,D)
  \leftarrow \Ext^n_{\Tri,\Ideal}(B,D)
  \leftarrow \Ext^n_{\Tri,\Ideal}(C,D)
  \leftarrow \Ext^{n-1}_{\Tri,\Ideal}(A,D) \leftarrow \\
  \dotsb
  \leftarrow \Ext^1_{\Tri,\Ideal}(C,D)
  \leftarrow \Ext^0_{\Tri,\Ideal}(A,D)
  \leftarrow \Ext^0_{\Tri,\Ideal}(B,D)
  \leftarrow \Ext^0_{\Tri,\Ideal}(C,D)
  \leftarrow 0
\end{multline*}
for any \(\Ideal\)\nb-exact exact triangle \(A\to B\to C\to \Sigma A\)
and any \(D\inOb\Tri\).

We claim that there are similar long exact sequences
\begin{multline*}
  0
  \to \Ext^0_{\Tri,\Ideal}(D,A)
  \to \Ext^0_{\Tri,\Ideal}(D,B)
  \to \Ext^0_{\Tri,\Ideal}(D,C)
  \to \Ext^1_{\Tri,\Ideal}(D,A)
  \to \dotsb\\
  \to \Ext^{n-1}_{\Tri,\Ideal}(D,C)
  \to \Ext^n_{\Tri,\Ideal}(D,A)
  \to \Ext^n_{\Tri,\Ideal}(D,B)
  \to \Ext^n_{\Tri,\Ideal}(D,C)
  \to \dotsb
\end{multline*}
in the second variable.  Since \(P(D)_n\) is \(\Ideal\)\nb-projective,
the sequences
\[
0 \to \Tri(P(D)_n,A) \to \Tri(P(D)_n,B) \to \Tri(P(D)_n,C) \to 0
\]
are exact for all \(n\in\N\).  This extension of chain complexes yields
the desired long exact sequence.

We list a few more elementary properties of derived functors.  We only
spell things out for the left derived functors \(\Left_n F\colon
\Tri\to\Cat\) of a covariant functor \(F\colon \Tri\to\Cat\).  Similar
assertions hold for right derived functors of contravariant functors.

The derived functors \(\Left_n F\) satisfy \(\Ideal\subseteq \ker
\Left_n F\) and hence descend to functors \(\Left_n F\colon
\Tri/\Ideal\to\Cat\) because the zero map \(P(A)\to P(B)\) is a chain
map lifting of~\(f\) if \(f\in\Ideal(A,B)\).  As a consequence,
\(\Left_n F(A)\cong0\) if~\(A\) is \(\Ideal\)\nb-contractible.  The
long exact homology sequences of Lemma~\ref{lem:derived_functors_les}
show that \(\Left_n F(f)\colon \Left_n F(A)\to \Left_n F(B)\) is
invertible if \(f\in\Tri(A,B)\) is an \(\Ideal\)\nb-equivalence.

\begin{warning}
  \label{warn:derived_not_homological}
  The derived functors \(\Left_n F\) are \emph{not homological} and
  therefore do not deserve to be called \(\Ideal\)\nb-exact even
  though they vanish on \(\Ideal\)\nb-phantom maps.
  Lemma~\ref{lem:derived_functors_les} shows that these functors are
  only half-exact on \(\Ideal\)\nb-exact triangles.  Thus \(\Left_n
  F(f)\) need not be monic (or epic) if~\(f\) is \(\Ideal\)\nb-monic
  (or \(\Ideal\)\nb-epic).  The problem is that the
  \(\Ideal\)\nb-projective resolution functor \(P\colon
  \Tri\to\Ho(\Tri)\) is not exact because it even fails to be stable.
\end{warning}

The following remarks require a more advanced background in homological
algebra and are not going to be used in the sequel.

\begin{remark}
  \label{rem:derived_functors_categories}
  The derived functors introduced above, especially the \(\Ext\)
  functors, can be interpreted in terms of \emph{derived categories}.

  We have already observed in~\S\ref{sec:homological_algebra_chain} that
  the \(\Ideal\)\nb-exact chain complexes form a thick subcategory of
  \(\Ho(\Tri)\).  The augmentation map \(P(A)\to A\) of an
  \(\Ideal\)\nb-projective resolution of \(A\inOb\Tri\) is a
  quasi-isomorphism with respect to this thick subcategory.  The chain
  complex \(P(A)\) is projective (see~\cite{Keller:Handbook}), that is,
  for any chain complex~\(C_\bullet\), the space of morphisms \(A\to
  C_\bullet\) in the derived category \(\Der(\Tri,\Ideal)\) agrees with
  \([P(A),C_\bullet]\).  Especially, \(\Ext^n_{\Tri,\Ideal}(A,B)\) is
  the space of morphisms \(A\to B[n]\) in \(\Der(\Tri,\Ideal)\).

  Now let \(F\colon \Tri\to\Cat\) be an additive covariant functor.
  Extend it to an exact functor \(\bar{F}\colon \Ho(\Tri)\to\Ho(\Cat)\).
  It has a total left derived functor
  \[
  \Left \bar{F}\colon \Der(\Tri,\Ideal)\to\Der(\Cat),
  \qquad A\mapsto \bar{F}\bigl(P(A)\bigr).
  \]
  By definition, we have \(\Left_n F(A) \defeq H_n\bigl(\Left
  \bar{F}(A)\bigr)\).
\end{remark}

\begin{remark}
  \label{rem:derived_cup_products}
  In classical Abelian categories, the \(\Ext\) groups form a graded
  ring, and the derived functors form graded modules over this graded
  ring.  The same happens in our context.  The most conceptual
  construction of these products uses the description of derived
  functors sketched in Remark~\ref{rem:derived_functors_categories}.

  Recall that we may view elements of \(\Ext^n_{\Tri,\Ideal}(A,B)\) as
  morphisms \(A\to B[n]\) in the derived category \(\Der(\Tri,\Ideal)\).
  Taking translations, we can also view them as morphisms \(A[m]\to
  B[n+m]\) for any \(m\in\Z\).  The usual composition in the category
  \(\Der(\Tri,\Ideal)\) therefore yields an associative product
  \[
  \Ext^n_{\Tri,\Ideal}(B,C)\otimes\Ext^m_{\Tri,\Ideal}(A,B)
  \to \Ext^{n+m}_{\Tri,\Ideal}(A,C).
  \]
  Thus we get a graded additive category with morphism spaces
  \(\bigl(\Ext^n_{\Tri,\Ideal}(A,B)\bigr)_{n\in\N}\).

  Similarly, if \(F\colon \Tri\to\Cat\) is an additive functor and
  \(\Left \bar{F}\colon \Der(\Tri,\Ideal)\to\Der(\Cat)\) is as in
  Remark~\ref{rem:derived_cup_products}, then a morphism \(A\to B[n]\)
  in \(\Der(\Tri,\Ideal)\) induces a morphism \(\Left\bar{F}(A)\to
  \Left\bar{F}(B)[n]\) in \(\Der(\Cat)\).  Passing to homology, we get
  canonical maps
  \[
  \Ext^n_{\Tri,\Ideal}(A,B) \to
  \Hom_\Cat\bigl(\Left F_m(A),\Left F_{m-n}(B)\bigr)
  \qquad \forall m\ge n,
  \]
  which satisfy an appropriate associativity condition.  For a
  contravariant functor, we get canonical maps
  \[
  \Ext^n_{\Tri,\Ideal}(A,B) \to
  \Hom_\Cat\bigl(\Right F^m(B),\Right F^{m+n}(A)\bigr)
  \qquad \forall m\ge 0.
  \]
\end{remark}

\subsection{Projective objects via adjointness}
\label{sec:projectives_adjointness}

We develop a method for constructing enough projective objects.  Let
\(\Tri\) and~\(\Cat\) be stable additive categories, let \(F\colon
\Tri\to\Cat\) be a stable additive functor, and let \(\Ideal\defeq \ker
F\).  In our applications, \(\Tri\) is triangulated and the
functor~\(F\) is either exact or stable and homological.

Recall that a covariant functor \(R\colon \Tri\to\Ab\) is
\emph{(co)representable} if it is naturally isomorphic to
\(\Tri(A,\blank)\) for some \(A\inOb\Tri\), which is then unique.  If
the functor \(B\mapsto \Cat\bigl(A,F(B)\bigr)\) on~\(\Tri\) is
representable, we write \(F^\lad(A)\) for the representing object.
By construction, we have natural isomorphisms
\[
\Tri(F^\lad(A),B)\cong \Cat\bigl(A,F(B)\bigr)
\]
for all \(B\in\Tri\).  Let~\(\Cat'\) be the full subcategory of all
objects \(A\inOb\Cat\) for which \(F^\lad(A)\) is defined.
Then~\(F^\lad\) is a functor \(\Cat'\to\Tri\), which we call the
\emph{(partially defined) left adjoint} of~\(F\).  Although one usually
assumes \(\Cat=\Cat'\), we shall also need~\(F^\lad\) in cases where it
is not defined everywhere.

The functor \(B\mapsto \Cat\bigl(A,F(B)\bigr)\) for \(A\inOb\Cat'\)
vanishes on \(\Ideal=\ker F\) for trivial reasons.  Hence
\(F^\lad(A)\inOb\Tri\) is \(\Ideal\)\nb-projective.  This simple
observation is surprisingly powerful: as we shall see, it often yields
all \(\Ideal\)\nb-projective objects.

\begin{remark}
  \label{rem:adjoint_properties}
  We have \(F^\lad(\Sigma A)\cong \Sigma F^\lad(A)\) for all
  \(A\inOb\Cat'\), so that \(\Sigma(\Cat')=\Cat'\).  Moreover,
  \(F^\lad\) commutes with infinite direct sums (as far as they exist
  in~\(\Tri\)) because
  \[
  \Tri\left( \bigoplus F^\lad(A_i), B\right) \cong \prod \Tri(
  F^\lad(A_i), B) \cong \prod \Cat\bigl(A_i, F(B)\bigr) \cong
  \Cat\left( \bigoplus A_i, F(B)\right).
  \]
\end{remark}

\begin{example}
  \label{exa:Ideal_K_projective}
  Consider the functor \(\K_*\colon \KK\to\Ab^{\Ztwo}\).  Let
  \(\Z\inOb\Ab^{\Ztwo}\) denote the trivially graded Abelian
  group~\(\Z\).  Notice that
  \begin{align*}
    \Hom\bigl(\Z,\K_*(A)\bigr) &\cong \K_0(A) \cong \KK(\C,A),\\
    \Hom\bigl(\Z[1],\K_*(A)\bigr) &\cong \K_1(A) \cong
    \KK(\CONT_0(\R),A),
  \end{align*}
  where \(\Z[1]\) means~\(\Z\) in odd degree.  Hence
  \(\K_*^\lad(\Z)=\C\) and \(\K_*^\lad(\Z[1])=\CONT_0(\R)\).  More
  generally, Remark~\ref{rem:adjoint_properties} shows that
  \(\K_*^\lad(A)\) is defined if both the even and odd parts of
  \(A\inOb\Ab^{\Ztwo}\) are countable free Abelian groups: it is a
  direct sum of at most countably many copies of \(\C\) and
  \(\CONT_0(\R)\).  Hence all such countable direct sums are
  \(\Ideal_\K\)\nb-projective (we briefly say
  \emph{\(\K\)\nb-projective}).  As we shall see, \(\K_*^\lad\) is
  not defined on all of \(\Ab^\Ztwo\); this is typical of homological
  functors.
\end{example}

\begin{example}
  \label{exa:ideal_H_projective}
  Consider the functor \(\Hgy\colon \Ho(\Cat;\Z/p)\to\Cat^{\Z/p}\) of
  Example~\ref{exa:homology}.  Let \(j\colon
  \Cat^{\Z/p}\to\Ho(\Cat;\Z/p)\) be the functor that views an object of
  \(\Cat^{\Z/p}\) as a \(p\)\nb-periodic chain complex whose boundary
  map vanishes.

  A chain map \(j(A)\to B_\bullet\) for \(A\inOb\Cat^{\Z/p}\) and
  \(B_\bullet\inOb\Ho(\Cat;\Z/p)\) is a family of maps \(\varphi_n\colon
  A_n\to \ker(d_n\colon B_n\to B_{n-1})\).  Such a family is chain
  homotopic to~\(0\) if and only if each~\(\varphi_n\) lifts to a map
  \(A_n\to B_{n+1}\).  Suppose that~\(A_n\) is projective for all
  \(n\in\Z/p\).  Then such a lifting exists if and only if
  \(\varphi_n(A_n)\subseteq d_{n+1}(B_{n+1})\).  Hence
  \[
  [j(A),B_\bullet] \cong
  \prod_{n\in\Z/p} \Cat\bigl(A_n,H_n(B_\bullet)\bigr) \cong
  \Cat^{\Z/p}\bigl(A,\Hgy(B_\bullet)\bigr).
  \]

  As a result, the left adjoint of~\(\Hgy\) is defined on the
  subcategory of projective objects \(\Proj(\Cat)^{\Z/p} \subseteq
  \Cat^{\Z/p}\) and agrees there with the restriction of~\(j\).  We will
  show in~\S\ref{sec:der_Ho_Cat} that \(\Proj(\Cat)^{\Z/p}\) is equal to
  the domain of~\(\Hgy^\lad\) and that all \(\IdealH\)\nb-projective
  objects are of the form \(\Hgy^\lad(A)\) for some
  \(A\inOb\Proj(\Cat)^{\Z/p}\) (provided~\(\Cat\) has enough projective
  objects).

  By duality, analogous results hold for injective objects: the domain
  of the \emph{right} adjoint of~\(\Hgy\) is the subcategory of injective
  objects of~\(\Cat^{\Z/p}\), the right adjoint is equal to~\(j\) on
  this subcategory, and this provides all \(\Hgy\)\nb-injective objects
  of \(\Ho(\Cat;\Z/p)\).
\end{example}

These examples show that~\(F^\lad\) yields many \(\ker
F\)-projective objects.  We want to get \emph{enough} \(\ker
F\)-projective objects in this fashion, assuming that~\(F^\lad\)
is defined on enough of~\(\Cat\).  In order to treat ideals of the form
\(\bigcap F_i\), we now consider a more complicated setup.  Let
\(\{\Cat_i\mid i\in I\}\) be a set of stable homological or triangulated
categories together with full subcategories
\(\Proj\Cat_i\subseteq\Cat_i\) and stable homological or exact functors
\(F_i\colon \Tri\to\Cat_i\) for all \(i\in I\).  Assume that
\begin{itemize}
\item the left adjoint~\(F_i^\lad\) is defined on~\(\Proj\Cat_i\) for
  all \(i\in I\);

\item there is an epimorphism \(P\to F_i(A)\) in~\(\Cat_i\) with
  \(P\inOb \Proj\Cat_i\) for any \(i\in I\), \(A\inOb\Tri\);

\item the set of functors \(F_i^\lad\colon \Proj\Cat_i\to\Tri\) is
  \emph{cointegrable}, that is, \(\bigoplus_{i\in I} F_i^\lad(B_i)\)
  exists for all families of objects \(B_i\in \Proj\Cat_i\), \(i\in I\).

\end{itemize}
The reason for the notation \(\Proj\Cat_i\) is that for a homological
functor~\(F_i\) we usually take~\(\Proj\Cat_i\) to be the class of
projective objects of~\(\Cat_i\); if~\(F_i\) is exact, then we often
take \(\Proj\Cat_i=\Cat_i\).  But it may be useful to choose a smaller
category, as long as it satisfies the second condition above.

\begin{proposition}
  \label{pro:epi_dagger}
  In this situation, there are enough \(\Ideal\)\nb-projective
  objects, and~\(\Proj_\Ideal\) is generated by \(\bigcup_{i\in I}
  \{F_i^\lad(B)\mid B\in \Proj\Cat_i\}\).  More precisely, an object
  of~\(\Tri\) is \(\Ideal\)\nb-projective if and only if it is a
  retract of \(\bigoplus_{i\in I} F_i^\lad(B_i)\) for a family of
  objects \(B_i\in \Proj\Cat_i\).
\end{proposition}

\begin{proof}
  Let \(\tilde\Proj_0 \defeq \bigcup_{i\in I} \{F_i^\lad(B)\mid B\in
  \Proj\Cat_i\}\) and \(\Proj_0\defeq \sumclosure{\tilde\Proj_0}\).  To
  begin with, we observe that any object of the form \(F_i^\lad(B)\)
  with \(B\inOb \Proj\Cat_i\) is \(\ker F_i\)\nb-projective and hence
  \(\Ideal\)\nb-projective because \(\Ideal\subseteq \ker F_i\).
  Hence~\(\Proj_0\) consists of \(\Ideal\)\nb-projective objects.

  Let \(A\inOb\Tri\).  For each \(i\in I\), there is an epimorphism
  \(p_i\colon B_i\to F_i(A)\) with \(B_i\in \Proj\Cat_i\).  The direct
  sum \(B\defeq \bigoplus_{i\in I} F_i^\lad(B_i)\) exists.  We have
  \(B\inOb\Proj_0\) by construction.  We are going to construct an
  \(\Ideal\)\nb-epimorphism \(p\colon B\to A\).  This shows that there
  are enough \(\Ideal\)\nb-projective objects.

  The maps \(p_i\colon B_i\to F_i(A)\) provide maps \(\hat{p}_i\colon
  F_i^\lad(B_i)\to A\) via the adjointness isomorphisms
  \(\Tri(F_i^\lad(B_i),A)\cong \Cat_i\bigl(B_i,F_i(A)\bigr)\).  We
  let \(p\defeq \sum \hat{p}_i\colon \bigoplus F_i^\lad(B_i)\to A\).
  We must check that~\(p\) is an \(\Ideal\)\nb-epimorphism.
  Equivalently, \(p\) is \(\ker F_i\)-epic for all \(i\in I\); this
  is, in turn equivalent to \(F_i(p)\) being an epimorphism
  in~\(\Cat_i\) for all \(i\in I\), because of Lemma \ref{lem:epic_kerF}
  or~\ref{lem:epic_kerF_exact}.  This is what we are going to prove.

  The identity map on \(F_i^\lad(B_i)\) yields a map \(\alpha_i\colon
  B_i\to F_iF_i^\lad(B_i)\) via the adjointness isomorphism
  \(\Tri\bigl(F_i^\lad(B_i),F_i^\lad(B_i)\bigr)\cong
  \Cat_i\bigl(B_i,F_iF_i^\lad(B_i)\bigr)\).  Composing with the map
  \[
  F_iF_i^\lad(B_i) \to
  F_i\left( \bigoplus F_i^\lad(B_i) \right) =
  F_i(B)
  \]
  induced by the coordinate embedding \(F_i^\lad(B_i)\to B\), we get
  a map \(\alpha_i'\colon B_i\to F_i(B)\).  The naturality of the
  adjointness isomorphisms yields \(F_i(\hat{p}_i)\circ\alpha_i=p_i\)
  and hence \(F_i(p)\circ \alpha_i'= p_i\).  The map~\(p_i\) is an
  epimorphism by assumption.  Now we use a cancellation result for
  epimorphisms: if \(f\circ g\) is an epimorphism, then so is~\(f\).
  Thus \(F_i(p)\) is an epimorphism as desired.

  If~\(A\) is \(\Ideal\)\nb-projective, then the
  \(\Ideal\)\nb-epimorphism \(p\colon B\to A\) splits; to see this,
  embed~\(p\) in an exact triangle \(N\to B\to A\to \Sigma N\) and
  observe that the map \(A\to\Sigma N\) belongs to \(\Ideal(A,\Sigma
  N)=0\).  Therefore, \(A\) is a retract of~\(B\).  Since~\(\Proj_0\) is
  closed under retracts and \(B\inOb\Proj_0\), we get \(A\inOb\Proj_0\).
  Hence~\(\tilde\Proj_0\) generates all \(\Ideal\)\nb-projective
  objects.
\end{proof}

\subsection{The universal exact homological functor}
\label{sec:universal_functor}

For the following results, it is essential to define an ideal by a
single functor~\(F\) instead of a family of functors as in
Proposition~\ref{pro:epi_dagger}.

\begin{definition}
  \label{def:universal_homological}
  Let~\(\Ideal\) be a homological ideal in a triangulated
  category~\(\Tri\).  An \(\Ideal\)\nb-exact stable homological functor
  \(F\colon \Tri\to\Cat\) is called \emph{universal} if any other
  \(\Ideal\)\nb-exact stable homological functor \(G\colon
  \Tri\to\Cat'\) factors as \(\bar{G}=G\circ F\) for a stable exact
  functor \(\bar{G}\colon \Cat\to\Cat'\) that is unique up to natural
  isomorphism.
\end{definition}

This universal property characterises~\(F\) uniquely up to natural
isomorphism.  We have constructed such a functor
in~\S\ref{sec:homological_algebra_chain}.  Beligiannis constructs it in
\cite{Beligiannis:Relative}*{\S3} using a localisation of the Abelian
category \(\Coh(\Tri)\) at a suitable Serre subcategory; he calls this
functor \emph{projectivisation functor} and its target category
\emph{Steenrod category}.  This notation is motivated by the special
case of the Adams spectral sequence.  The following theorem allows us to
check whether a given functor is universal:

\begin{theorem}
  \label{the:universal_functor}
  Let~\(\Tri\) be a triangulated category, let \(\Ideal\subseteq\Tri\)
  be a homological ideal, and let \(F\colon \Tri\to\Cat\) be an
  \(\Ideal\)\nb-exact stable homological functor into a stable Abelian
  category~\(\Cat\); let~\(\Proj\Cat\) be the class of projective
  objects in~\(\Cat\).  Suppose that idempotent morphisms in~\(\Tri\)
  split.

  The functor~\(F\) is the universal \(\Ideal\)\nb-exact stable
  homological functor and there are enough \(\Ideal\)\nb-projective
  objects in~\(\Tri\) if and only if
  \begin{itemize}
  \item \(\Cat\) has enough projective objects;

  \item the adjoint functor~\(F^\lad\) is defined on~\(\Proj\Cat\);

  \item \(F\circ F^\lad(A)\cong A\) for all \(A\inOb\Proj\Cat\).

  \end{itemize}
\end{theorem}

\begin{proof}
  Suppose first that~\(F\) is universal and that there are enough
  \(\Ideal\)\nb-projective objects.  Then~\(F\) is equivalent to the
  projectivisation functor of~\cite{Beligiannis:Relative}.  The various
  properties of this functor listed in
  \cite{Beligiannis:Relative}*{Proposition 4.19} include the following:
  \begin{itemize}
  \item there are enough projective objects in~\(\Cat\);

  \item \(F\) induces an equivalence of categories
    \(\Proj_\Ideal\cong\Proj\Cat\) (\(\Proj_\Ideal\) is the class of
    projective objects in~\(\Tri\));

  \item \(\Cat\bigl(F(A),F(B)\bigr) \cong \Tri(A,B)\) for all
    \(A\inOb\Proj_\Ideal\), \(B\inOb\Tri\).

  \end{itemize}
  Here we use the assumption that idempotents in~\(\Tri\) split.  The
  last property is equivalent to \(F^\lad\circ F(A)\cong A\) for all
  \(A\inOb\Proj_\Ideal\).  Since \(\Proj_\Ideal\cong\Proj\Cat\)
  via~\(F\), this implies that~\(F^\lad\) is defined on \(\Proj\Cat\)
  and that \(F\circ F^\lad(A)\cong A\) for all \(A\inOb\Proj\Cat\).
  Thus~\(F\) has the properties listed in the statement of the theorem.

  Now suppose conversely that~\(F\) has these properties.  Let
  \(\Proj_\Ideal'\subseteq\Tri\) be the essential range of
  \(F^\lad\colon \Proj\Cat\to \Tri\).  We claim that~\(\Proj_\Ideal'\)
  is the class of all \(\Ideal\)\nb-projective objects in~\(\Tri\).
  Since \(F\circ F^\lad\) is equivalent to the identity functor
  on~\(\Proj\Cat\) by assumption, \(F|_{\Proj_\Ideal'}\) and~\(F^\lad\)
  provide an equivalence of categories \(\Proj_\Ideal'\cong \Proj\Cat\).
  Since~\(\Cat\) is assumed to have enough projectives, the hypotheses
  of Proposition~\ref{pro:epi_dagger} are satisfied.  Hence there are
  enough \(\Ideal\)\nb-projective objects in~\(\Tri\), and any object
  of~\(\Proj_\Ideal\) is a retract of an object of~\(\Proj_\Ideal'\).
  Idempotent morphisms in the category \(\Proj_\Ideal'\cong\Proj\Cat\)
  split because~\(\Cat\) is Abelian and retracts of projective objects
  are again projective.  Hence~\(\Proj_\Ideal'\) is closed under
  retracts in~\(\Tri\), so that \(\Proj_\Ideal'=\Proj_\Ideal\).  It also
  follows that \(F\) and~\(F^\lad\) provide an equivalence of categories
  \(\Proj_\Ideal\cong\Proj\Cat\).  Hence \(F^\lad\circ F(A)\cong A\) for
  all \(A\inOb\Proj_\Ideal\), so that we get \(\Cat\bigl(F(A),F(B)\bigr)
  \cong \Tri(F^\lad\circ F(A),B) \cong \Tri(A,B)\) for all
  \(A\inOb\Proj_\Ideal\), \(B\inOb\Tri\).

  Now let \(G\colon \Tri\to\Cat'\) be a stable homological functor.  We
  will later assume~\(G\) to be \(\Ideal\)\nb-exact, but the first
  part of the following argument works in general.  Since~\(F\) provides
  an equivalence of categories \(\Proj_\Ideal\cong\Proj\Cat\), the rule
  \(\bar{G}\bigl(F(P)\bigr)\defeq G(P)\) defines a functor~\(\bar{G}\)
  on \(\Proj\Cat\).  This yields a functor \(\Ho(\bar{G})\colon
  \Ho(\Proj\Cat)\to\Ho(\Cat')\).  Since~\(\Cat\) has enough projective
  objects, the construction of projective resolutions provides a functor
  \(P\colon \Cat\to \Ho(\Proj\Cat)\).  We let~\(\bar{G}\) be the
  composite functor
  \[
  \bar{G}\colon \Cat \xrightarrow{P}
  \Ho(\Proj\Cat) \xrightarrow{\Ho(\bar{G})}
  \Ho(\Cat') \xrightarrow{H_0}
  \Cat'.
  \]
  This functor is right-exact and satisfies \(\bar{G}\circ F=G\) on
  \(\Ideal\)\nb-projective objects of~\(\Tri\).

  Now suppose that~\(G\) is \(\Ideal\)\nb-exact.  Then we get
  \(\bar{G}\circ F=G\) for all objects of~\(\Tri\) because this holds
  for \(\Ideal\)\nb-projective objects.  We claim that~\(\bar{G}\) is
  exact.  Let \(A\inOb\Cat\).  Since~\(\Cat\) has enough projective
  objects, we can find a projective resolution of~\(A\).  We may assume
  this resolution to have the form \(F(P_\bullet)\) with
  \(P_\bullet\inOb\Ho(\Proj_\Ideal)\) because \(F(\Proj_\Ideal)\cong
  \Proj\Cat\).  Lemma~\ref{lem:exact_kerF} yields that~\(P_\bullet\) is
  \(\Ideal\)\nb-exact except in degree~\(0\).  Since \(\Ideal
  \subseteq \ker G\), the chain complex~\(P_\bullet\) is \(\ker
  G\)-exact in positive degrees as well, so that \(G(P_\bullet)\)
  is exact except in degree~\(0\) by Lemma~\ref{lem:exact_kerF}.  As a
  consequence, \(\Left_p \bar{G}(A)=0\) for all \(p>0\).  We also have
  \(\Left_0\bar{G}(A)=\bar{G}(A)\) by construction.  Thus~\(\bar{G}\) is
  exact.

  As a result, \(G\) factors as \(G=\bar{G}\circ F\) for an exact
  functor \(\bar{G}\colon \Cat\to\Cat'\).  It is clear that~\(\bar{G}\)
  is stable.  Finally, since~\(\Cat\) has enough projective objects, a
  functor on~\(\Cat\) is determined up to natural equivalence by its
  restriction to projective objects.  Therefore, our factorisation
  of~\(G\) is unique up to natural equivalence.  Thus~\(F\) is the
  universal \(\Ideal\)\nb-exact functor.
\end{proof}

\begin{remark}
  \label{rem:adjoint_defined_on_some_projectives}
  Let \(\Proj'\Cat\subseteq\Proj\Cat\) be some subcategory such that any
  object of~\(\Cat\) is a quotient of a direct sum of objects
  of~\(\Proj'\Cat\).  Equivalently,
  \(\sumclosure{\Proj'\Cat}=\Proj\Cat\).
  Theorem~\ref{the:universal_functor} remains valid if we only assume
  that~\(F^\lad\) is defined on~\(\Proj'\Cat\) and that \(F\circ
  F^\lad(A)\cong A\) holds for \(A\inOb\Proj'\Cat\) because both
  conditions are evidently hereditary for direct sums and retracts.
\end{remark}

\begin{theorem}
  \label{the:universal_homological_nice}
  In the situation of Theorem~\textup{\ref{the:universal_functor}},
  the domain of the functor~\(F^\lad\) is equal to~\(\Proj\Cat\), and
  its essential range is~\(\Proj_\Ideal\).  The functors \(F\)
  and~\(F^\lad\) restrict to equivalences of categories
  \(\Proj_\Ideal\cong\Proj\Cat\) inverse to each other.

  An object \(A\inOb\Tri\) is \(\Ideal\)\nb-projective if and only if
  \(F(A)\) is projective and
  \[
  \Cat\bigl(F(A),F(B)\bigr) \cong \Tri(A,B)
  \]
  for all \(B\inOb\Tri\); following Ross
  Street~\cite{Street:Homotopy_classification}, we call such objects
  \emph{\(F\)\nb-projective}.  We have \(F(A)\inOb\Proj\Cat\) if and
  only if there is an \(\Ideal\)\nb-equivalence \(P\to A\) with
  \(P\inOb\Proj_\Ideal\).

  The functors \(F\) and~\(F^\lad\) induce bijections between
  isomorphism classes of projective resolutions of~\(F(A)\) in~\(\Cat\)
  and isomorphism classes of \(\Ideal\)\nb-projective resolutions of
  \(A\inOb\Tri\) in~\(\Tri\).

  If \(G\colon \Tri\to\Cat'\) is any \textup{(}stable\textup{)}
  homological functor, then there is a unique right-exact
  \textup{(}stable\textup{)} functor \(\bar{G}\colon \Cat\to\Cat'\) such
  that \(\bar{G}\circ F(P)=G(P)\) for all \(P\inOb\Proj_\Ideal\).

  The left derived functors of~\(G\) with respect to~\(\Ideal\) and
  of~\(\bar{G}\) are related by natural isomorphisms
  \(\Left_n\bar{G}\circ F(A)=\Left_n G(A)\) for all \(A\inOb\Tri\),
  \(n\in\N\).  There is a similar statement for cohomological functors,
  which specialises to natural isomorphisms
  \[
  \Ext^n_{\Tri,\Ideal}(A,B) \cong \Ext^n_\Cat\bigl(F(A),F(B)\bigr).
  \]
\end{theorem}

\begin{proof}
  We have already seen during the proof of
  Theorem~\ref{the:universal_functor} that~\(F\) restricts to an
  equivalence of categories \(\Proj_\Ideal\congto\Proj\Cat\) with
  inverse~\(F^\lad\) and that \(\Cat\bigl(F(A),F(B)\bigr)\cong
  \Tri(A,B)\) for all \(A\inOb\Proj_\Ideal\), \(B\inOb\Proj\Cat\).

  Conversely, if~\(A\) is \(F\)\nb-projective in the sense of Street,
  then~\(A\) is \(\Ideal\)\nb-projective because already
  \(\Tri(A,B)\cong \Cat\bigl(F(A),F(B)\bigr)\) for all \(B\inOb\Tri\)
  yields \(A\cong F^\lad\circ F(A)\), so that~\(A\) is
  \(\Ideal\)\nb-projective; notice that the projectivity of \(F(A)\)
  is automatic.

  Since~\(F\) maps \(\Ideal\)\nb-equivalences to isomorphisms,
  \(F(A)\) is projective whenever there is an
  \(\Ideal\)\nb-equivalence \(P\to A\) with
  \(\Ideal\)\nb-projective~\(P\).  Conversely, suppose that \(F(A)\)
  is \(\Ideal\)\nb-projective.  Let \(P_0\to A\) be a one-step
  \(\Ideal\)\nb-projective resolution.  Since \(F(A)\) is projective,
  the epimorphism \(F(P_0)\to F(A)\) splits by some map \(F(A)\to
  F(P_0)\).  The resulting map \(F(P_0)\to F(A)\to F(P_0)\) is
  idempotent and comes from an idempotent endomorphism of~\(P_0\)
  because~\(F\) is fully faithful on~\(\Proj_\Ideal\).  Its range
  object~\(P\) exists because we require idempotent morphisms
  in~\(\Cat\) to split.  It belongs again to~\(\Proj_\Ideal\), and the
  induced map \(F(P)\to F(A)\) is invertible by construction.  Hence we
  get an \(\Ideal\)\nb-equivalence \(P\to A\).

  If~\(C_\bullet\) is a chain complex over~\(\Tri\), then we know
  already from Lemma~\ref{lem:exact_kerF} that~\(C_\bullet\) is
  \(\Ideal\)\nb-exact if and only if \(F(C_\bullet)\) is exact.
  Hence~\(F\) maps an \(\Ideal\)\nb-projective resolution of~\(A\) to
  a projective resolution of \(F(A)\).  Conversely, if \(P_\bullet\to
  F(A)\) is any projective resolution in~\(\Cat\), then it is of the
  form \(F(\hat{P}_\bullet)\to F(A)\) where \(\hat{P}_\bullet\defeq
  F^\lad(P_\bullet)\) and where we get the map \(\hat{P}_0\to A\) by
  adjointness from the given map \(P_0\to F(A)\).  This shows that~\(F\)
  induces a bijection between isomorphism classes of
  \(\Ideal\)\nb-projective resolutions of~\(A\) and projective
  resolutions of \(F(A)\).

  We have seen during the proof of Theorem~\ref{the:universal_functor}
  how a stable homological functor \(G\colon \Tri\to\Cat'\) gives rise
  to a unique right-exact functor \(\bar{G}\colon \Cat\to\Cat'\) that
  satisfies \(\bar{G}\circ F(P)=G(P)\) for all \(P\inOb\Proj_\Ideal\).
  The derived functors \(\Left_n\bar{G}\bigl(F(A)\bigr)\) for
  \(A\inOb\Tri\) are computed by applying~\(\bar{G}\) to a projective
  resolution of \(F(A)\).  Since such a projective resolution is of the
  form \(F(P_\bullet)\) for an \(\Ideal\)\nb-projective resolution
  \(P_\bullet\to A\) and since \(\bar{G}\circ F=G\) on
  \(\Ideal\)\nb-projective objects, the derived functors \(\Left_n
  G(A)\) and \(\Left_n\bar{G}\bigl(F(A)\bigr)\) are computed by the same
  chain complex and agree.  The same reasoning applies to cohomological
  functors and yields the assertion about \(\Ext\).

  Finally, we check that \(A\inOb\Cat\) is projective if
  \(F^\lad(A)\) is defined.  We prove \(\Ext^1_\Cat(A,B)=0\) for all
  \(B\inOb\Cat\), from which the assertion follows.  There is a
  projective resolution of the form \(F(P_\bullet)\to B\), which we use
  to compute \(\Ext^1_\Cat(A,B)\).  The adjointness of~\(F^\lad\)
  and~\(F\) yields that \(F^\lad(A)\inOb\Tri\) is
  \(\Ideal\)\nb-projective and that \(\Cat\bigl(A,
  F(P_\bullet)\bigr)\cong \Tri(F^\lad(A),P_\bullet)\).
  Since~\(P_\bullet\) is \(\Ideal\)\nb-exact in positive degrees by
  Lemma~\ref{lem:exact_kerF} and~\(F^\lad(A)\) is
  \(\Ideal\)\nb-projective, we get \(0 =
  H_1\bigl(\Cat\bigl(A,F(P_\bullet)\bigr)\bigr)=\Ext^1_\Cat(A,B)\).
\end{proof}

\begin{remark}
  \label{rem:idempotents_split}
  The assumption that idempotents split is only needed to check that the
  universal \(\Ideal\)\nb-exact functor has the properties listed in
  Theorem~\ref{the:universal_functor}.  The converse directions of
  Theorem~\ref{the:universal_functor} and
  Theorem~\ref{the:universal_homological_nice} do not need this
  assumption.

  If~\(\Tri\) has countable direct sums or countable direct products,
  then idempotents in~\(\Tri\) automatically split by
  \cite{Neeman:Triangulated}*{\S1.3}.  This covers categories
  such as \(\KK^G\) because they have countable direct sums.
\end{remark}

\subsection{Derived functors in homological algebra}
\label{sec:der_Ho_Cat}

Now we study the kernel~\(\IdealH\) of the homology functor \(\Hgy\colon
\Ho(\Cat;\Z/p)\to \Cat^{\Z/p}\) introduced Example~\ref{exa:homology}.
We get exactly the same statements if we replace the homotopy category
by its derived category and study the kernel of \(\Hgy\colon
\Der(\Cat;\Z/p)\to \Cat^{\Z/p}\).  We often abbreviate~\(\IdealH\)
to~\(\HGY\) and speak of \(\HGY\)\nb-epimorphisms, \(\HGY\)\nb-exact
chain complexes, \(\HGY\)\nb-projective resolutions, and so on.  We
denote the full subcategory of \(\HGY\)\nb-projective objects in
\(\Ho(\Cat;\Z/p)\) by~\(\Proj_{\HGY}\).

We assume that the underlying Abelian category~\(\Cat\) has enough
projective objects.  Then the same holds for~\(\Cat^{\Z/p}\), and we
have \(\Proj(\Cat^{\Z/p}) \cong (\Proj\Cat)^{\Z/p}\).  That is, an
object of \(\Cat^{\Z/p}\) is projective if and only if its homogeneous
pieces are.

\begin{theorem}
  \label{the:homology_universal}
  The category \(\Ho(\Cat;\Z/p)\) has enough \(\HGY\)\nb-projective
  objects, and the functor \(\Hgy\colon \Ho(\Cat;\Z/p)\to \Cat^{\Z/p}\)
  is the universal \(\HGY\)\nb-exact stable homological functor.  Its
  restriction to~\(\Proj_{\HGY}\) provides an equivalence of categories
  \(\Proj_{\HGY}\cong \Proj\Cat^{\Z/p}\).  More concretely, a chain
  complex in \(\Ho(\Cat;\Z/p)\) is \(\HGY\)\nb-projective if and only
  if it is homotopy equivalent to one with vanishing boundary map and
  projective entries.

  The functor~\(\Hgy\) maps isomorphism classes of
  \(\HGY\)\nb-projective resolutions of a chain complex~\(A\) in
  \(\Ho(\Cat;\Z/p)\) bijectively to isomorphism classes of projective
  resolutions of \(\Hgy(A)\) in~\(\Cat^{\Z/p}\).  We have
  \[
  \Ext^n_{\Ho(\Cat;\Z/p),\IdealH}(A,B) \cong
  \Ext^n_\Cat\bigl(\Hgy(A),\Hgy(B)\bigr).
  \]
  Let \(F\colon \Cat\to\Cat'\) be some covariant additive functor and
  define
  \[
  \bar{F}\colon \Ho(\Cat;\Z/p)\to \Ho(\Cat';\Z/p)
  \]
  by applying~\(F\) entrywise.  Then \(\Left_n\bar{F}(A) \cong \Left_n
  F\bigl(\Hgy(A)\bigr)\) for all \(n\in\N\).  Similarly, we have
  \(\Right^n\bar{F}(A) \cong \Right^n F\bigl(\Hgy(A)\bigr)\) if~\(F\) is
  a contravariant functor.
\end{theorem}

\begin{proof}
  The category~\(\Cat^{\Z/p}\) has enough projective objects by
  assumption.  We have already seen in
  Example~\ref{exa:ideal_H_projective} that~\(\Hgy^\lad\) is defined
  on \(\Proj\Cat^{\Z/p}\); this functor is denoted by~\(j\) in
  Example~\ref{exa:ideal_H_projective}.  It is clear that \(\Hgy\circ
  j(A)\cong A\) for all \(A\inOb\Cat^{\Z/p}\).  Now
  Theorem~\ref{the:universal_functor} shows that~\(\Hgy\) is universal.
  We do not need idempotent morphisms in \(\Ho(\Cat;\Z/p)\) to split by
  Remark~\ref{rem:idempotents_split}.
\end{proof}

\begin{remark}
  \label{rem:ideal_versus_truncation}
  Since the universal \(\Ideal\)\nb-exact functor is essentially
  unique, the universality of \(\Hgy\colon
  \Der(\Cat;\Z/p)\to\Cat^{\Z/p}\) means that we can recover this functor
  and hence the stable Abelian category~\(\Cat^{\Z/p}\) from the ideal
  \(\IdealH\subseteq \Der(\Cat;\Z/p)\).  That is, the
  ideal~\(\IdealH\) and the functor \(\Hgy\colon
  \Der(\Cat;\Z/p)\to\Cat^{\Z/p}\) contain exactly the same amount of
  information.

  For instance, if we forget the precise category~\(\Cat\) by
  composing~\(\Hgy\) with some \emph{faithful} functor \(\Cat\to\Cat'\),
  then the resulting homology functor \(\Ho(\Cat;\Z/p)\to\Cat'\) still
  has kernel~\(\IdealH\).  We can recover~\(\Cat^{\Z/p}\) by passing to
  the universal \(\Ideal\)\nb-exact functor.

  We compare this with the situation for truncation structures
  (\cite{Beilinson-Bernstein-Deligne}).  These cannot exist for periodic
  categories such as \(\Der(\Cat;\Z/p)\) for \(p\ge1\).  Given the
  standard truncation structure on \(\Der(\Cat)\), we can recover the
  Abelian category~\(\Cat\) as its core; we also get back the homology
  functors \(H_n\colon \Der(\Cat)\to\Cat\) for all \(n\in\Z\).
  Conversely, the functor \(\Hgy\colon \Der(\Cat)\to\Cat^\Z\) together
  with the grading on~\(\Cat^\Z\) tells us what it means for a chain
  complex to be exact in degrees greater than~\(0\) or less than~\(0\)
  and thus determines the truncation structure.  Hence the standard
  truncation structure on \(\Der(\Cat)\) contains the same amount of
  information as the functor \(\Hgy\colon \Der(\Cat)\to\Cat^\Z\)
  together with the grading on~\(\Cat^\Z\).
\end{remark}

\section{The plain Universal Coefficient Theorem}
\label{sec:plain_UCT}

Now we study the ideal \(\Ideal_\K\defeq \ker \K_*\subseteq\KK\) of
Example~\ref{exa:Kvan}.  We complete our analysis of this example and
explain the Universal Coefficient Theorem for~\(\KK\) in our framework.
We call \(\Ideal_\K\)\nb-projective objects and \(\Ideal_\K\)\nb-exact
functors briefly \emph{\(\K\)\nb-projective} and \emph{\(\K\)\nb-exact}
and let~\(\Proj_\K\) be the class of \(\K\)\nb-projective objects
in~\(\KK\).

Let \(\CAb^\Ztwo\subseteq \Ab^\Ztwo\) be the full subcategory of
\emph{countable} \(\Ztwo\)\nb-graded Abelian groups.  Since the
\(\K\)\nb-theory of a separable \(\Cst\)\nb-algebra is countable, we
may view~\(\K_*\) as a stable homological functor \(\K_*\colon \KK\to
\CAb^\Ztwo\).

\begin{theorem}
  \label{the:K_projectives}
  There are enough \(\K\)\nb-projective objects in \(\KK\), and the
  universal \(\K\)\nb-exact functor is \(\K_*\colon \KK\to
  \CAb^\Ztwo\).  It restricts to an equivalence of categories between
  \(\Proj_\K\)~and the full subcategory \(\FCAb^\Ztwo\subseteq
  \CAb^\Ztwo\) of \(\Ztwo\)\nb-graded countable free Abelian groups.
  A separable \(\Cst\)\nb-algebra belongs to~\(\Proj_\K\) if and only
  if it is \(\KK\)-equivalent to \(\bigoplus_{i\in I_0} \C \oplus
  \bigoplus_{i\in I_1} \CONT_0(\R)\) where the sets \(I_0,I_1\) are at
  most countable.

  If \(A\inOb\KK\), then~\(\K_*\) maps isomorphism classes of
  \(\K\)\nb-projective resolutions of~\(A\) in~\(\Tri\) bijectively to
  isomorphism classes of free resolutions of \(\K_*(A)\).  We have
  \[
  \Ext^n_{\KK,\Ideal_\K}(A,B) \cong
  \begin{cases}
    \Hom_{\Ab^{\Ztwo}}\bigl(\K_*(A),\K_*(B)\bigr) & \text{for \(n=0\);}\\
    \Ext^1_{\Ab^{\Ztwo}}\bigl(\K_*(A),\K_*(B)\bigr) & \text{for \(n=1\);}\\
    0 & \text{for \(n\ge2\).}\\
  \end{cases}
  \]

  Let \(F\colon \KK\to\Cat\) be some covariant additive functor; then
  there is a unique right-exact functor \(\bar{F}\colon
  \CAb^\Ztwo\to\Cat\) with \(\bar{F}\circ\K_*=F\).  We have \(\Left_n F=
  (\Left_n\bar{F})\circ \K_*\) for all \(n\in\N\); this vanishes for
  \(n\ge2\).  Similar assertions hold for contravariant functors.
\end{theorem}

\begin{proof}
  Notice that \(\CAb^\Ztwo\subseteq \Ab^\Ztwo\) is an Abelian category.
  We shall denote objects of \(\Ab^\Ztwo\) by pairs \((A_0,A_1)\) of
  Abelian groups.  By definition, \((A_0,A_1)\inOb \FCAb^\Ztwo\) if and
  only if \(A_0\) and~\(A_1\) are countable free Abelian groups, that
  is, they are of the form \(A_0 = \Z[I_0]\) and \(A_1 = \Z[I_1]\) for
  at most countable sets \(I_0, I_1\).  It is well-known that any
  Abelian group is a quotient of a free Abelian group and that subgroups
  of free Abelian groups are again free.  Moreover, free Abelian groups
  are projective.  Hence \(\FCAb^\Ztwo\) is the subcategory of
  projective objects in \(\CAb^\Ztwo\) and any object \(G\inOb
  \CAb^\Ztwo\) has a projective resolution of the form \(0 \to F_1\to
  F_0 \prto G\) with \(F_0,F_1\inOb \FCAb^\Ztwo\).  This implies that
  derived functors on \(\CAb^\Ztwo\) only occur in dimensions \(1\)
  and~\(0\).

  As in Example~\ref{exa:Ideal_K_projective}, we see
  that~\(\K_*^\lad\) is defined on \(\FCAb^\Ztwo\) and satisfies
  \[
  \K_*^\lad\bigl(\Z[I_0],\Z[I_1]\bigr) \cong
  \bigoplus_{i\in I_0} \C \oplus \bigoplus_{i\in I_1} \CONT_0(\R)
  \]
  if \(I_0,I_1\) are countable.  We also have \(\K_*\circ
  \K_*^\lad\bigl(\Z[I_0],\Z[I_1]\bigr) \cong
  \bigl(\Z[I_0],\Z[I_1]\bigr)\), so that the hypotheses of
  Theorem~\ref{the:universal_functor} are satisfied.  Hence there are
  enough \(\K\)\nb-projective objects and~\(\K_*\) is universal.  The
  remaining assertions follow from
  Theorem~\ref{the:universal_homological_nice} and our detailed
  knowledge of the homological algebra in \(\CAb^\Ztwo\).
\end{proof}

\begin{example}
  \label{exa:Kuenneth_derived}
  Consider the stable homological functor
  \[
  F\colon \KK\to \CAb^\Ztwo,
  \qquad
  A\mapsto \K_*(A\otimes B)
  \]
  for some \(B\inOb\KK\), where~\(\otimes\) denotes, say, the spatial
  \(\Cst\)\nb-tensor product.  We claim that the associated right-exact
  functor \(\CAb^\Ztwo\to \CAb^\Ztwo\) is
  \[
  \bar{F}\colon \CAb^\Ztwo\to \CAb^\Ztwo,
  \qquad
  G \mapsto G\otimes \K_*(B).
  \]
  It is easy to check \(F\circ \K_*^\lad(G) \cong G\otimes\K_*(B)
  \cong \bar{F}(G)\) for \(G\inOb \FCAb^\Ztwo\).  Since the functor
  \(G\mapsto G\otimes \K_*(B)\) is right-exact and agrees
  with~\(\bar{F}\) on projective objects, we get \(\bar{F}(G) =
  G\otimes\K_*(B)\) for all \(G\inOb \CAb^\Ztwo\).  Hence the derived
  functors of~\(F\) are
  \[
  \Left_n F(A) \cong
  \begin{cases}
    \K_*(A)\otimes\K_*(B) &\text{for \(n=0\);}\\
    \Tor^1\bigl(\K_*(A),\K_*(B)\bigr) &\text{for \(n=1\);}\\
    0  &\text{for \(n\ge2\).}\\
  \end{cases}
  \]
  Here we use the same graded version of \(\Tor\) as in the K\"unneth
  Theorem (\cite{Blackadar:K-theory}).
\end{example}

\begin{example}
  \label{exa:KK_second_derived}
  Consider the stable homological functor
  \[
  F\colon \KK\to \Ab^\Ztwo,
  \qquad B\mapsto \KK_*(A,B)
  \]
  for some \(A\inOb\KK\).  We suppose that~\(A\) is a \emph{compact}
  object of~\(\KK\), that is, the functor~\(F\) commutes with direct
  sums.  Then \(\KK_*\bigl(A,\K_*^\lad(G)\bigr) \cong
  \KK_*(A,\C)\otimes G\) for all \(G\inOb \FCAb^\Ztwo\) because this
  holds for \(G=(\Z,0)\) and is inherited by suspensions and direct
  sums.  Now we get \(\bar{F}(G) \cong \KK_*(A,\C) \otimes G\) for all
  \(G\inOb \CAb^\Ztwo\) as in Example~\ref{exa:Kuenneth_derived}.
  Therefore,
  \[
  \Left_n F(B) \cong
  \begin{cases}
    \KK_*(A,\C)\otimes\K_*(B) &\text{for \(n=0\);}\\
    \Tor^1\bigl(\KK_*(A,\C),\K_*(B)\bigr) &\text{for \(n=1\);}\\
    0  &\text{for \(n\ge2\).}\\
  \end{cases}
  \]
\end{example}

Generalising Examples \ref{exa:Kuenneth_derived}
and~\ref{exa:KK_second_derived}, we have \(\bar{F}(G) \cong F(\C)\otimes
G\) and hence
\[
\Left_n F(B) \cong
\begin{cases}
  F(\C)\otimes\K_*(B) &\text{for \(n=0\),}\\
  \Tor^1\bigl(F(\C),\K_*(B)\bigr) &\text{for \(n=1\),}
\end{cases}
\]
for any covariant functor \(F\colon \KK\to\Cat\) that commutes with
direct sums.

Similarly, if \(F\colon \KK\to\Cat\) is contravariant and maps direct
sums to direct products, then \(\bar{F}(G) \cong \Hom(G,F(\C))\) and
\[
\Right^n F(B) \cong
\begin{cases}
  \Hom\bigl(\K_*(B), F(\C)\bigr) &\text{for \(n=0\),}\\
  \Ext^1\bigl(\K_*(B), F(\C)\bigr) &\text{for \(n=1\).}
\end{cases}
\]
The description of \(\Ext^n_{\KK,\Ideal_\K}\) in
Theorem~\ref{the:K_projectives} is a special case of this.

\subsection*{Universal Coefficient Theorem in the hereditary case}
%\addcontentsline{toc}{subsection}{Universal Coefficient Theorem in the hereditary case}
\label{sec:UCT_1dim}

In general, we need spectral sequences in order to relate the derived
functors~\(\Left_n F\) back to~\(F\).  We will discuss this in a sequel
to this article.  Here we only treat the simple case where we have
projective resolutions of length~\(1\).  The following universal
coefficient theorem is very similar to but slightly more general than
\cite{Beligiannis:Relative}*{Theorem 4.27} because we do not require
\emph{all} \(\Ideal\)\nb-equivalences to be invertible.

\begin{theorem}
  \label{the:UCT_homological}
  Let \(\Ideal\) be a homological ideal in a triangulated
  category~\(\Tri\).  Let \(A\inOb\Tri\) have an
  \(\Ideal\)\nb-projective resolution of length~\(1\).  Suppose also
  that \(\Tri(A,B)=0\) for all \(\Ideal\)\nb-contractible~\(B\).  Let
  \(F\colon \Tri\to\Cat\) be a homological functor, \(\tilde{F}\colon
  \Tri^\op\to\Cat\) a cohomological functor, and \(B\inOb\Tri\).  Then
  there are natural short exact sequences
  \begin{gather*}
    0 \to \Left_0 F_*(A) \to F_*(A) \to \Left_1 F_{*-1}(A) \to 0,\\
    0 \to \Right^1 \tilde{F}^{*-1}(A) \to \tilde{F}^*(A) \to
    \Right^0 \tilde{F}^*(A) \to 0,\\
    0 \to \Ext^1_{\Tri,\Ideal}(\Sigma A,B) \to \Tri(A,B) \to
    \Ext^0_{\Tri,\Ideal}(A,B) \to 0.
  \end{gather*}
\end{theorem}

\begin{example}
  \label{exa:KK_UCT}
  For the ideal \(\Ideal_\K\subseteq\KK\), any object has a
  \(\K\)\nb-projective resolution of length~\(1\) by
  Theorem~\ref{the:K_projectives}.  The other hypothesis of
  Theorem~\ref{the:UCT_homological} holds if and only if~\(A\) satisfies
  the Universal Coefficient Theorem (UCT).  The UCT for \(\KK(A,B)\)
  predicts \(\KK(A,B)=0\) if \(\K_*(B)=0\).  Conversely, if this is the
  case, then Theorem~\ref{the:UCT_homological} applies, and our
  description of \(\Ext_{\KK,\Ideal_\K}\) in
  Theorem~\ref{the:K_projectives} yields the UCT for \(\KK(A,B)\) for
  all~\(B\).  This yields our claim.

  Thus the UCT for \(\KK(A,B)\) is a special of
  Theorem~\ref{the:UCT_homological}.  In the situations of Examples
  \ref{exa:Kuenneth_derived} and~\ref{exa:KK_second_derived}, we get the
  familiar K\"unneth Theorems for \(\K_*(A\otimes B)\) and
  \(\KK_*(A,B)\).  These arguments are very similar to the original
  proofs (see~\cite{Blackadar:K-theory}).  Our machinery allows us to
  treat other situations in a similar fashion.
\end{example}

\begin{proof}[Proof of Theorem~\ref{the:UCT_homological}]
  We only write down the proof for homological functors.  The
  cohomological case is dual and contains \(\Tri(\blank,B)\) as a
  special case.

  Let \(0\to P_1\xrightarrow{\delta_1} P_0\xrightarrow{\delta_0} A\) be
  an \(\Ideal\)\nb-projective resolution of length~\(1\) and view it
  as an \(\Ideal\)\nb-exact chain complex of length~\(3\).
  Lemma~\ref{lem:exact_complex_length_3} yields a commuting diagram
  \[\xymatrix{
    P_1 \ar[r]^{\delta_1} \ar@{=}[d] &
    P_0 \ar@{=}[d] \ar[r]^{\tilde{\delta}_0} &
    \tilde{A} \ar@{.>}[d]^{\alpha} \\
    P_1 \ar[r]^{\delta_1} &
    P_0 \ar[r]^{\delta_0} &
    A,
  }
  \]
  such that the top row is part of an \(\Ideal\)\nb-exact exact
  triangle \(P_1\to P_0 \to \tilde{A} \to \Sigma P_1\) and~\(\alpha\) is
  an \(\Ideal\)\nb-equivalence.  We claim that~\(\alpha\) is an
  isomorphism in~\(\Tri\).

  We embed~\(\alpha\) in an exact triangle \(\Sigma^{-1}B \to
  \tilde{A} \xrightarrow{\alpha} A\xrightarrow{\beta} B\).
  Lemma~\ref{lem:equivalence_contractible} shows that \(B\) is
  \(\Ideal\)\nb-contractible because~\(\alpha\) is an
  \(\Ideal\)\nb-equivalence.  Hence \(\Tri(A,B)=0\) by our
  assumption on~\(A\).  This forces \(\beta=0\), so that our exact
  triangle splits: \(\tilde{A}\cong A \oplus \Sigma^{-1}B\).  Now we
  apply the functor $\Tri(\cdot,B)$ to the exact triangle \(P_0\to
  P_1\to\tilde{A}\).  The resulting long exact sequence has the form
  \[
  \dotsb \leftarrow \Tri(P_0,B) \leftarrow \Tri(\tilde{A},B)
  \leftarrow \Tri(\Sigma P_1,B) \leftarrow \dotsb.
  \]
  Since both \(P_0\) and~\(P_1\) are \(\Ideal\)\nb-projective
  and~\(B\) is \(\Ideal\)\nb-contractible, we get
  \(\Tri(\tilde{A},B)=0\).  Then
  \(\Tri(B,B)\subseteq\Tri(\tilde{A},B)\) vanishes as well, so that
  \(B\cong 0\) and~\(\alpha\) is invertible.

  We get an exact triangle in~\(\Tri\) of the form \(P_1
  \xrightarrow{\delta_1} P_0 \xrightarrow{\delta_0} A \to \Sigma P_1\)
  because any triangle isomorphic to an exact one is itself exact.

  Now we apply~\(F\).  Since~\(F\) is homological, we get a long exact
  sequence
  \[
  \dotsb \to
  F_*(P_1) \xrightarrow{F_*(\delta_1)}
  F_*(P_0) \to
  F_*(A) \to
  F_{*-1}(P_1) \xrightarrow{F_{*-1}(\delta_1)}
  F_{*-1}(P_0) \to \dotsb.
  \]
  We cut this into short exact sequences of the form
  \[
  \coker \bigl(F_*(\delta_1)\bigr) \into
  F_*(A) \prto
  \ker \bigl(F_{*-1}(\delta_1)\bigr).
  \]
  Since \(\coker F_*(\delta_1)=\Left_0 F_*(A)\) and \(\ker
  F_*(\delta_1)=\Left_1 F_*(A)\), we get the desired exact sequence.
  The map \(\Left_0 F_*(A)\to F_*(A)\) is the canonical map induced
  by~\(\delta_0\).  The other map \(F_*(A)\to \Left_1 F_{*-1}(A)\) is
  natural for all morphisms between objects with an
  \(\Ideal\)\nb-projective resolution of length~\(1\) by
  Proposition~\ref{pro:resolutions}.
\end{proof}

The proof shows that -- in the situation of
Theorem~\ref{the:UCT_homological} -- we have
\[
\Ext^0_{\Tri,\Ideal}(A,B) \cong \Tri/\Ideal(A,B),
\qquad
\Ext^1_{\Tri,\Ideal}(A,B) \cong \Ideal(A,\Sigma B).
\]
More generally, we can construct a natural map \(\Ideal(A,\Sigma B)\to
\Ext^1_{\Tri,\Ideal}(A,B)\) for any homological ideal, using the
\(\Ideal\)\nb-universal homological functor \(F\colon \Tri\to\Cat\).
We embed \(f\in\Ideal(A,\Sigma B)\) in an exact triangle \(B\to C\to
A\to \Sigma B\).  We get an extension
\[
\bigl[F(B) \into F(C) \prto F(A)\bigr] \inOb
\Ext^1_\Cat\bigl(F(A), F(B)\bigr)
\]
because this triangle is \(\Ideal\)\nb-exact.  This class
\(\kappa(f)\) in \(\Ext^1_\Cat\bigl(F(A), F(B)\bigr)\) does not depend
on auxiliary choices because the exact triangle \(B\to C\to A\to \Sigma
B\) is unique up to isomorphism.
Theorem~\ref{the:universal_homological_nice} yields
\(\Ext^1_{\Tri,\Ideal}(A, B) \cong \Ext^1_\Cat\bigl(F(A), F(B)\bigr)\)
because~\(F\) is universal.  Hence we get a natural map
\[
\kappa\colon \Ideal(A,\Sigma B) \to \Ext^1_{\Tri,\Ideal}(A,B).
\]
We may view~\(\kappa\) as a secondary invariant generated by the
canonical map
\[
\Tri(A,B)\to\Ext^0_{\Tri,\Ideal}(A,B).
\]
For the ideal~\(\Ideal_\K\), we get the same map~\(\kappa\) as in
Example~\ref{exa:Kvan}.

An Abelian category with enough projective objects is called
\emph{hereditary} if any subobject of a projective object is again
projective.  Equivalently, any object has a projective resolution of
length~\(1\).  This motivates the following definition:

\begin{definition}
  \label{def:hereditary}
  A homological ideal~\(\Ideal\) in a triangulated category~\(\Tri\) is
  called \emph{hereditary} if any object of~\(\Tri\) has a projective
  resolution of length~\(1\).
\end{definition}

If~\(\Ideal\) is hereditary and if \(\Ideal\)\nb-equivalences are
invertible, then Theorem~\ref{the:UCT_homological} applies to all
\(A\inOb\Tri\) (and vice versa).

\section{Crossed products for compact quantum groups}
\label{sec:crossed_cqg}

Let~\(G\) be a compact (quantum) group and write \(\CONT(G)\) and
\(\Cst(G)\) for the Hopf \(\Cst\)\nb-algebras of \(\hat{G}\) and~\(G\).
We study the homological algebra in~\(\KK^G\) generated by the ideals
\(\Ideal_{\ltimes,\K}\subseteq \Ideal_{\ltimes}\subseteq \KK\) defined
in Example~\ref{exa:cross}.  Recall that~\(\Ideal_{\ltimes}\) is the
kernel of the crossed product functor
\[
G\cross\blank\colon \KK^G\to \KK,
\qquad A\mapsto G\cross A,
\]
whereas \(\Ideal_{\ltimes,\K}\) is the kernel of the composite functor
\(\K_*\circ (G\cross\blank)\).  Since we have already analysed~\(\K_*\)
in~\S\ref{sec:plain_UCT}, we can treat both ideals in a parallel
fashion.

Our setup contains two classical special cases.  First, \(G\) may be a
compact Lie group.  Then \(\CONT(G)\) and \(\Cst(G)\) have the usual
meaning, and the objects of~\(\KK^G\) are separable
\(\Cst\)\nb-algebras with a continuous action of~\(G\).  Secondly,
\(\CONT(G)\) may be the dual \(\Cst_\red(H)\) of a discrete group~\(H\),
so that \(\Cst(G)=\CONT_0(H)\).  Then the objects of~\(\KK^G\) are
separable \(\Cst\)\nb-algebras equipped with a (reduced) coaction
of~\(H\).  (We disregard the nuances between reduced and full
coactions.)  If we identify \(\KK^{\Cst_\red(H)} \cong \KK^H\) using
Baaj-Skandalis duality, then the crossed product functor
\(G\cross\blank\) corresponds to the forgetful functor \(\KK^H\to\KK\).

Let \(\Proj_\ltimes\) and \(\Proj_{\ltimes,\K}\) be the classes of
projective objects for the ideals \(\Ideal_\ltimes\)
and~\(\Ideal_{\ltimes,\K}\).  Our first task is to find enough
projective objects for these two ideals.

Let \(\tau\colon \KK\to\KK^G\) be the functor that equips \(B\inOb\KK\)
with the trivial \(G\)\nb-action (that is, coaction of \(\CONT(G)\)).
This corresponds to the induction functor from the trivial quantum group
to \(\Cst(G)\) under the equivalence \(\KK^G\cong \KK^{\Cst(G)}\).

The functors \(\tau\) and \(G\cross\blank\) are adjoint, that is, there
are natural isomorphisms
\begin{equation}
  \label{eq:tau_cross_adjoint}
  \KK^G(\tau(A),B)\cong \KK(A,G\cross B)
\end{equation}
for all \(A\inOb\KK\), \(B\inOb\KK^G\).  This generalisation of the
Green-Julg Theorem is proved in \cite{Vergnioux:These}*{Th\'eor\`eme 5.10}.

\begin{lemma}
  \label{lem:Ideal_cross_projectives}
  There are enough projective objects for \(\Ideal_\ltimes\) and
  \(\Ideal_{\ltimes,\K}\).  We have
  \[
  \Proj_\ltimes = \sumclosure[\bigl]{\tau(\KK)},
  \qquad
  \Proj_{\ltimes,\K} =
  \sumclosure[\bigl]{\tau\C, \tau \CONT_0(\R)}.
  \]
\end{lemma}

\begin{proof}
  The adjoint of \(G\cross\blank\) is defined on all of~\(\KK\), which
  is certainly enough for Proposition~\ref{pro:epi_dagger}.  This yields
  the assertions for \(\Ideal_\ltimes\) and~\(\Proj_\ltimes\); even
  more, \(\Proj_\ltimes\) is the closure of \(\tau(\KK)\) under
  retracts.

  The adjoint of \(\K_*\circ (G\cross\blank)\) is \(\tau\circ
  \K_*^\lad\), which is defined for free countable
  \(\Ztwo\)\nb-graded Abelian groups.  Explicitly,
  \[
  \tau\circ \K_*^\lad\bigl(\Z[I_0],\Z[I_1]\bigr)
  \cong \bigoplus_{i\in I_0} \tau(\C) \oplus \bigoplus_{i\in I_1}
  \tau\bigl(\CONT_0(\R)\bigr).
  \]
  Since any object of \(\CAb^\Ztwo\) is a quotient of one in
  \(\FCAb^\Ztwo\), Proposition~\ref{pro:epi_dagger} applies and yields
  the assertions about \(\Ideal_{\ltimes,\K}\)
  and~\(\Proj_{\ltimes,\K}\); even more, \(\Proj_{\ltimes,\K}\) is the
  closure of \(\tau\circ \K_*^\lad(\FCAb^\Ztwo)\) under retracts.
\end{proof}

To proceed further, we describe the universal \(\Ideal\)\nb-exact
functors.  The functors \((G\cross \blank)\) and \(\K_*\circ
(G\cross\blank)\) fail the criterion of
Theorem~\ref{the:universal_functor} (unless \(G=\{1\}\)) because
\[
G\cross \tau(A) \cong \Cst(G) \otimes A \not\cong A.
\]
This is not surprising because \(G\cross\blank\) is equivalent to a
\emph{forgetful} functor.  The universal functor recovers a
\emph{linearisation} of the forgotten structure.

First we consider the ideal~\(\Ideal_{\ltimes,\K}\).  By
Lemma~\ref{lem:Ideal_cross_projectives}, the objects \(\tau(\C)\) and
\(\Sigma\tau(\C)\) generate all \(\Ideal_{\ltimes,\K}\)-projective
objects.  Their internal symmetries are encoded by the
\(\Ztwo\)\nb-graded ring \(\KK^G_*(\tau\C,\tau\C)^\op\); the
superscript~\(\op\) denotes that we take the product in reversed order.
For classical compact groups, this ring is commutative, so that the
order of multiplication does not matter; in general, the reversed-order
product is the more standard choice.  The following fact is well-known:

\begin{lemma}
  \label{lem:Rep_via_KK}
  The ring \(\KK^G_*(\tau\C,\tau\C)^\op\) is isomorphic to the
  \emph{representation ring} \(\Rep(G)\) of~\(G\) for \(*=0\) and~\(0\)
  for \(*=1\).
\end{lemma}

We may take this as the definition of \(\Rep(G)\).

\begin{proof}
  The adjointness isomorphism~\eqref{eq:tau_cross_adjoint} identifies
  \[
  \KK^G_*(\tau\C,\tau\C) \cong \KK_*(\C,G\cross \tau\C) = \KK_*\bigl(\C,
  \Cst(G)\bigr) \cong \K_*(\Cst G).
  \]
  Since~\(G\) is compact, \(\Cst G\) is a direct sum of matrix algebras,
  one for each irreducible representation of~\(G\).  Hence the
  underlying Abelian group of \(\Rep(G)\) is the free Abelian group
  \(\Z[\hat{G}]\) on the set~\(\hat{G}\) of irreducible representations
  of~\(G\).  The ring structure on \(\Rep(G)\) comes from the internal
  tensor product of representations: represent two elements
  \(\alpha,\beta\) of \(\KK^G_0(\tau\C,\tau\C)\) by differences of
  finite-dimensional representations \(\pi,\varrho\) of~\(G\), then
  \(\alpha\circ\beta\in\KK^G_0(\tau\C,\tau\C)\) is represented by
  \(\varrho\otimes\pi\) because the product in~\(\KK^G\) boils down to
  an exterior tensor product in this case (with reversed order).
\end{proof}

\begin{example}
  \label{exa:Rep_dual_group}
  If \(\CONT(G)=\Cst_\red(H)\) for a discrete group~\(H\), then
  \(\hat{G}=H\) and the product on \(\Rep(G)=\Z[H]\) is the usual
  convolution.  Thus \(\Rep(G)\) is the group ring of~\(H\).

  If~\(G\) is a compact group, then \(\Rep(G)\) is the
  representation ring in the usual sense.
\end{example}

For any \(B\inOb\KK^G\), the Kasparov product turns
\(\KK^G_*(\tau\C,B)\) into a left module over the ring
\(\KK^G_0(\tau\C,\tau\C)^\op\cong \Rep(G)\).  Thus \(\KK^G_*(\tau\C,B)\)
becomes an object of the Abelian category \(\CMod{\Rep G}^\Ztwo\) of
\(\Ztwo\)\nb-graded countable \(\Rep(G)\)-modules.  We get a
stable homological functor
\[
F_\K\colon \KK^G\to \CMod{\Rep G}^\Ztwo,
\qquad F_\K(B) \defeq \KK^G_*(\tau\C,B).
\]
By~\eqref{eq:tau_cross_adjoint}, the underlying Abelian group of
\(\KK^G_*(\tau\C,B)\) is
\[
\KK^G_*(\tau\C,B) \cong \KK_*(\C,G\cross B) \cong \K_*(G\cross B).
\]
Hence we still have \(\ker F_\K=\Ideal_{\ltimes,\K}\).  We often write
\(F_\K(B)= \K_*(G\cross B)\) if it is clear from the context that we
view \(\K_*(G\cross B)\) as an object of \(\CMod{\Rep G}^\Ztwo\).

\begin{theorem}
  \label{the:cross_K_universal}
  The functor~\(F_\K\) is the universal
  \(\Ideal_{\ltimes,\K}\)-exact functor.

  The subcategory of \(\Ideal_{\ltimes,\K}\)\nb-projective objects in
  \(\KK^G\) is equivalent to the subcategory of \(\Ztwo\)\nb-graded
  countable projective \(\Rep(G)\)-modules.  If \(A\inOb\KK^G\),
  then~\(F_\K\) induces a bijection between isomorphism classes of
  \(\Ideal_{\ltimes,\K}\)-projective resolutions of~\(A\) and
  projective resolutions of \(F_\K(A)\) in \(\CMod{\Rep G}^\Ztwo\).

  If \(A,B\inOb\KK^G\), then
  \[
  \Ext^n_{\KK^G,\Ideal_{\ltimes,\K}} (A,B) \cong
  \Ext^n_{\Rep(G)}\bigl(\K_*(G\cross A), \K_*(G\cross B)\bigr).
  \]

  If \(H\colon \KK^G\to\Cat\) is homological and commutes with countable
  direct sums, then
  \[
  \Left_n H(A) \cong
  \Tor_n^{\Rep(G)}\bigl(H_*(\tau\C),\K_*(G\cross A)\bigr);
  \]
  if \(H\colon (\KK^G)^\op\to\Cat\) is cohomological and turns countable
  direct sums into direct products, then
  \[
  \Right^n H(A) \cong
  \Ext^n_{\Rep(G)}\bigl(\K_*(G\cross A), H^*(\tau\C)\bigr);
  \]
  Here we use the right or left \(\Rep(G)\)-module structure on
  \(H_*(\tau\C)\) that comes from the functoriality of~\(H\).
\end{theorem}

\begin{proof}
  We verify universality using Theorem~\ref{the:universal_functor}.  The
  category \(\CMod{\Rep G}^\Ztwo\) has enough projective objects:
  countable free modules are projective, and any object is a quotient of
  a countable free module.

  Given a free module \((\Rep(G)[I_0],\Rep(G)[I_1])\), we have natural
  isomorphisms
  \begin{multline*}
    \Hom_{\Rep(G)}\bigl((\Rep(G)[I_0],\Rep(G)[I_1]), F_\K(B)\bigr)
    \cong \prod_{\varepsilon\in\{0,1\}, i\in I_\varepsilon}
    \K_\varepsilon(G\cross B)
    \\ \cong \KK^G\biggl(
    \bigoplus_{\varepsilon\in\{0,1\}, i\in I_\varepsilon}
    \Sigma^\varepsilon\tau\C, B\biggr).
  \end{multline*}
  Hence the adjoint functor~\(F_\K^\lad\) is defined on countable free
  modules.  Idempotents in~\(\KK^G\) split by
  Remark~\ref{rem:idempotents_split}.  Therefore, the domain
  of~\(F_\K^\lad\) is closed under retracts and contains all projective
  objects of \(\CMod{\Rep G}^\Ztwo\).  It is easy to see that
  \(F_\K\circ F_\K^\lad(A)\cong A\) for free modules.  This extends to
  retracts and hence holds for all projective modules (compare
  Remark~\ref{rem:adjoint_defined_on_some_projectives}).  Now
  Theorem~\ref{the:universal_functor} yields that~\(F_\K\) is universal.

  The assertions about projective objects, projective resolutions, and
  \(\Ext\) now follow from Theorem~\ref{the:universal_homological_nice}.
  Theorem~\ref{the:universal_homological_nice} also yields a formula for
  left derived functors in terms of the right-exact functor
  \(\bar{H}\colon \CMod{\Rep G}^\Ztwo\to\Cat\) associated to a
  homological functor \(H\colon \KK^G\to\Cat\).  It remains to
  compute~\(\bar{H}\).

  First we define the \(\Tor\) objects in the statement of the theorem
  if~\(\Cat\) is the category of Abelian groups.  Then
  \(H_*(\tau\C)\inOb \Mod{\Rep G}^\Ztwo\), and we can take the derived
  functors of the usual \(\Ztwo\)\nb-graded balanced tensor product
  \(\otimes_{\Rep G}\) for \(\Rep(G)\)-modules.  We claim that
  there are natural isomorphisms
  \[
  \bar{H}_*(M) \cong H_*(\tau\C)\otimes_{\Rep G} M
  \]
  for all \(M\inOb \CMod{\Rep G}^\Ztwo\).  This holds for \(M=\Rep G\)
  and hence for all free modules because we have natural isomorphisms
  \[
  \bar{H}_*(\Rep G)
  \cong H_*\bigl(F_\K^\lad(\Rep G)\bigr)
  \cong H_*(\tau\C)
  \cong H_*(\tau\C) \otimes_{\Rep G} \Rep G.
  \]
  For general~\(M\), the functor~\(\bar{H}(M)\) is computed by a free
  resolution because it is right-exact.  Using this, one extends the
  computation to all modules.  By definition, \(\Tor^n_{\Rep G}(N,M)\)
  for \(N\inOb \Mod{(\Rep G)^\op}^\Ztwo\), \(M\inOb \CMod{\Rep
  G}^\Ztwo\), is the \(n\)th left derived functor of the functor
  \(N\otimes_{\Rep G}\blank\) on \(\CMod{\Rep G}^\Ztwo\).

  Now Theorem~\ref{the:universal_homological_nice} yields the formula
  for \(\Left_n H\) provided~\(H\) takes values in Abelian groups.  The
  same argument works in general, we only need more complicated
  categories.

  Let \(\Cat^\Ztwo[\Rep(G)^\op]\) be the category of
  \(\Ztwo\)\nb-graded objects~\(A\) of~\(\Cat\) together with a ring
  homomorphism \(\Rep(G)^\op\to \Cat^\Ztwo(A,A)\).  We can extend the
  definition of \(\otimes_{\Rep G}\) to get an additive stable bifunctor
  \[
  \otimes_{\Rep G}\colon \Cat^\Ztwo[\Rep(G)^\op] \otimes \Mod{\Rep G}^\Ztwo
  \to \Cat^\Ztwo.
  \]
  Its derived functors are \(\Tor^n_{\Rep G}\).  As above, we see that
  this yields the derived functors of~\(H\).

  The case of cohomological functors is similar and left to the reader.
\end{proof}

If~\(G\) is a compact group, then the same derived functors appear in
the Universal Coefficient Theorem for~\(\KK^G\) by Jonathan Rosenberg
and Claude Schochet (\cite{Rosenberg-Schochet:Kunneth}).  This is no
coincidence, of course.  We postpone a further discussion because we do
not treat spectral sequences here.

In order to get the universal functor for the ideal~\(\Ideal_\ltimes\),
we must lift the \(\Rep G\)-module structure on \(\K_*(G\cross B)\) to
\(G\cross B\).  Given any additive category~\(\Cat\), we define a
category \(\Cat[\Rep G]\) as in the proof of
Theorem~\ref{the:cross_K_universal}: its objects are pairs~\((A,\mu)\)
where \(A\inOb\Cat\) and~\(\mu\) is a ring homomorphism
\(\Rep(G)\to\Cat(A,A)\); its morphisms are morphisms in~\(\Cat\) that
are compatible with the \(\Rep(G)\)-module structure in the obvious
sense.

A module structure \(\mu\colon \Rep(G)\to \KK(A,A)\) for \(A\inOb\KK\)
is equivalent to a natural family of \(\Rep(G)\)-module structures
in the usual sense on the groups \(\KK(D,A)\) for all \(D\inOb\KK\):
simply define \(x\cdot y\defeq \mu(x)\circ y\) for \(x\in\Rep(G)\),
\(y\in\KK(D,A)\) and notice that this recovers the homomorphism~\(\mu\)
for \(y=\ID_A\).

The crucial property of the universal \(\Ideal_{\ltimes,\K}\)\nb-exact
functor~\(F_\K\) is that it lifts the original functor
\(\K_*(G\cross\blank)\colon \KK^G\to \CAb^\Ztwo\) to a functor
\[
F_\K\colon \KK^G\to \CMod{\Rep G}^\Ztwo = \CAb^\Ztwo[\Rep G].
\]
We need a similar lifting of \(G\cross\blank\colon \KK^G\to\KK\) to
\(\KK[\Rep G]\).  This requires a simple special case of exterior
products in~\(\KK^G\).  In general, it is not so easy to define exterior
products in~\(\KK^G\) for quantum groups because diagonal actions on
\(\Cst\)\nb-algebras are not defined without additional structure.  The
only case where this is easy is if one of the factors carries the
\emph{trivial} coaction.  This exterior product operation on the level
of \(\Cst\)\nb-algebras also works for Kasparov cycles, that is, we get
canonical maps
\[
\KK^G_0(A,B) \otimes \KK_0(C,D) \to
\KK^G_0(A\otimes C, B\otimes D)
\]
for all \(A,B\inOb\KK^G\), \(C,D\inOb\KK\).  Equivalently,
\((A,C)\mapsto A\otimes C\) is a bifunctor \(\KK^G\otimes\KK\to\KK^G\).

This exterior product construction yields a natural map
\[
\varrho_A\colon \Rep(G)^\op\cong \KK^G(\tau\C,\tau\C) \to
\KK^G(\tau\C\otimes A,\tau\C\otimes A) \cong
\KK^G(\tau A,\tau A)
\]
for \(A\inOb\KK\), whose range commutes with the range of the map
\[
\tau\colon \KK(A,A)\to \KK^G(\tau A,\tau A).
\]

The ring \(\KK^G(\tau A,\tau A)\) acts on \(\KK^G(\tau A,B)\cong
\KK(A,G\cross B)\) on the right by Kasparov product.  Hence so does
\(\Rep(G)^\op\) via~\(\varrho_A\).  Thus \(\KK(A,G\cross B)\) becomes a
left \(\Rep(G)\)-module for all \(A\inOb\KK\), \(B\inOb\KK^G\).
These module structures are natural in the variable~\(A\) because the
images of \(\KK^G(\tau\C,\tau\C)\) and \(\KK(A,A)\) in \(\KK^G(\tau
A,\tau A)\) commute.  Hence they must come from a ring homomorphism
\[
\mu_B\colon \Rep(G)\to \KK(G\cross B,G\cross B).
\]
These ring homomorphisms are natural because the \(\Rep(G)\)-module
structures on \(\KK(A, G\cross B)\) are manifestly natural in~\(B\).
Thus we have lifted \(G\cross\blank\) to a functor
\[
F\colon \KK^G\to \KK[\Rep G],
\qquad
B\mapsto (G\cross B,\mu_B).
\]
It is clear that \(\ker F=\Ideal_\ltimes\).  The target category
\(\KK[\Rep G]\) is neither triangulated nor Abelian.  To remedy this, we
use the Yoneda embedding \(\Yoneda\colon \KK\to \Coh(\KK)\) constructed
in~\S\ref{sec:universal_homological}.  This embedding is fully faithful;
so is the resulting functor \(\KK[\Rep G]\to \Coh(\KK)[\Rep G]\).

\begin{theorem}
  \label{the:cross_universal}
  The functor \(\Yoneda\circ F\colon \KK^G\to \Coh(\KK)[\Rep G]\) is the
  universal \(\Ideal_\ltimes\)\nb-exact functor.
\end{theorem}

\begin{proof}
  We omit the proof of this theorem because it is only notationally more
  difficult than the proof of Theorem~\ref{the:cross_K_universal}.
\end{proof}

The category \(\Coh(\KK)[\Rep G]\) is not as terrible as it seems.  We
can usually stay within the more tractable subcategory \(\KK[\Rep G]\),
and many standard techniques of homological algebra like bar resolutions
work in this setting.  This often allows us to compute derived functors
on \(\Coh(\KK)[\Rep G]\) in more classical terms.

Recall that the bar resolution of a \(\Rep G\)-module~\(M\) is a
natural free resolution
\[
\dotsb \to (\Rep G)^{\otimes n} \otimes M
\to (\Rep G)^{\otimes n-1} \otimes M
\to \dotsb
\to \Rep G\otimes M
\to M
\]
with certain natural boundary maps.  Defining
\[
(\Rep G)^{\otimes n} \otimes M =
\Z[\hat{G}^n] \otimes M \defeq \bigoplus_{x\in \hat{G}^n} M,
\]
we can make sense of this in \(\Cat[\Rep G]\) provided~\(\Cat\) has
countable direct sums; the \(\Rep G\)-module structures and the
boundary maps can also be defined.

If \(A\inOb\KK[\Rep G]\), then the bar resolution lies in \(\KK[\Rep
G]\) and is a projective resolution of~\(A\) in the ambient Abelian
category \(\Coh(\KK)[\Rep G]\).  Hence we can use it to compute derived
functors.  For the extension groups, we get
\[
\Ext^n_{\KK^G,\Ideal_\ltimes}(A,B) \cong
\HH^n\bigl(\Rep(G); \KK(G\cross A,G\cross B) \bigr);
\]
here \(\HH^n(R;M)\) denotes the \(n\)th Hochschild cohomology of a
ring~\(R\) with coefficients in an \(R\)\nb-bimodule~\(M\), and
\(\KK(G\cross A,G\cross B)\) is a bimodule over \(\Rep G\) via the
Kasparov product on the left and right and the ring homomorphisms
\[
\Rep(G)\to \KK(G\cross A,G\cross A),
\qquad
\Rep(G)\to \KK(G\cross B,G\cross B).
\]
Similarly, if \(H\colon \KK^G\to\Ab\) commutes with direct sums, then
\[
\Left_n H(B) \cong \HH_n\bigl(\Rep(G); H\circ \tau (G\cross B)\bigr),
\]
where \(H\circ\tau (G\cross B)\) carries the following \(\Rep
G\)-bimodule structure: the left module structure comes from
\(\mu_B\colon \Rep(G)\to \KK(G\cross B,G\cross B)\) and the right one
comes from \(\varrho_A\colon \Rep(G)^\op\to \KK^G(\tau A,\tau A)\) for
\(A=G\cross B\) and the functoriality of~\(H\).  The details are left to
the reader.

\subsection*{The Pimsner-Voiculescu exact sequence}
%\addcontentsline{toc}{subsection}{The Pimsner-Voiculescu exact sequence}
\label{sec:PV}

The reader may wonder why we have considered the
ideal~\(\Ideal_\ltimes\), given that the derived functors
for~\(\Ideal_{\ltimes,\K}\) are so much easier to describe.  This is
related to the question whether \(\Ideal\)\nb-equivalences are
invertible.  The ideal~\(\Ideal_{\ltimes,\K}\) cannot have this property
because it already fails for trivial~\(G\).  In contrast, the
ideal~\(\Ideal_\ltimes\) sometimes has this property.  This means that
the spectral sequences that we get from~\(\Ideal_\ltimes\) may converge
for all objects, not just for those in an appropriate bootstrap
category.  To illustrate this, we explain how the well-known
Pimsner-Voiculescu exact sequence fits into our framework.

This exact sequence deals with actions of the group~\(\Z\); to remain in
the framework of~\S\ref{sec:crossed_cqg}, we use Baaj-Skandalis duality
to turn such actions into actions of the Pontrjagin dual
group~\(\Torus\).  The representation ring of~\(\Torus\) is the ring
\(R\defeq \Z[t,t^{-1}]\) of Laurent polynomials or, equivalently, the
group ring of~\(\Z\).  An \(R\)\nb-bimodule in an Abelian
category~\(\Cat\) is an object~\(M\) of~\(\Cat\) together with two
commuting automorphisms \(\lambda,\rho\colon M\to M\).  The Hochschild
homology and cohomology are easy to compute using the free bimodule
resolution \(0\to R^{\otimes 2} \overset{d}\to R^{\otimes 2}\to R\),
where \(d(f) = t\cdot f\cdot t^{-1}-f\).  We get
\begin{align*}
  \HH_0(R;M) &\cong \HH^1(R;M) \cong \coker (\lambda\rho^{-1}-1),\\
  \HH_1(R;M) &\cong \HH^0(R;M) \cong \ker (\lambda\rho^{-1}-1),
\end{align*}
and \(\HH_n(R;M)\cong \HH^n(R;M)\cong 0\) for \(n\ge2\).  Transporting
this kind of resolution to \(\Cat[R]\), we get that any object of
\(\Cat[R]\) has an \(\Ideal_\ltimes\)\nb-projective resolution of
length~\(1\).  This would fail for \(\Ideal_{\ltimes,\K}\) because the
category of \(R\)\nb-modules has a non-trivial \(\Ext^2\).

The crucial point is that \(\Ideal_\ltimes\)\nb-equivalences are
invertible in~\(\KK^{\Torus}\).  By Baaj-Skandalis duality, this is
equivalent to the following statement: if \(f\inOb\KK^\Z(A,B)\) becomes
invertible in~\(\KK\), then it is already invertible in~\(\KK^\Z\).  We
do not want to discuss here how to prove this.  Taking this for granted,
we can now apply Theorem~\ref{the:UCT_homological} to all objects of
\(\KK^{\Torus}\).

We write down the resulting exact sequences for \(\KK^\Z(A,B)\) for
\(A,B\inOb\KK^\Z\) because this equivalent setting is more familiar.
The actions of~\(\Z\) on \(A\) and~\(B\) provide two actions of~\(\Z\)
on \(\KK_*(A,B)\).  We let \(t_A, t_B\colon \KK_*(A,B)\to\KK_*(A,B)\) be
the actions of the generators.  Theorem~\ref{the:UCT_homological} yields
an exact sequence
\[
\coker\bigl(t_At_B^{-1}-1|_{\KK_{*+1}(A,B)}\bigr)
\into \KK^\Z_*(A,B)
\prto \ker\bigl(t_At_B^{-1}-1|_{\KK_*(A,B)}\bigr).
\]
This is equivalent to a long exact sequence
\[\xymatrix{
  \KK_1(A,B) \ar[r] &
  \KK_0^\Z(A,B) \ar[r] &
  \KK_0(A,B) \ar[d]^{t_At_B^{-1}-1} \\
  \KK_1(A,B) \ar[u]^{t_At_B^{-1}-1} &
  \KK_1^\Z(A,B) \ar[l] &
  \KK_0(A,B). \ar[l]
}
\]
Similar manipulations yield the Pimsner-Voiculescu exact sequence for
the functor \(A\mapsto \K_*(\Z\cross A)\) and more general functors
defined on~\(\KK^\Z\).

\end{document}